# Relations among Ramanujan-Type Congruences II

## Ramanujan-Type Congruences in Half-Integral Weights

Martin Raum*

**Abstract:** We link Ramanujan-type congruences, which emerge abundantly in combinatorics, to the Galois- and geometric theory of modular forms. Specifically, we show that Ramanujan-type congruences are preserved by the action of the shallow Hecke algebra, and discover a dichotomy between congruences originating in Hecke eigenvalues and congruences on arithmetic progressions with cube-free periods. The latter provide congruences among algebraic parts of twisted central L-values. We specialize our results to integer partitions, for which we investigate the landmark proofs of partition congruences by Atkin and by Ono. Based on a modulo $\ell$ analogue of the Maeda conjecture for certain partition generating functions, we conclude that their approach by Hecke operators acting diagonally modulo $\ell$ on modular forms is indeed close to optimal. This work is enabled by several structure results for Ramanujan-type congruences that we establish. In an extended example, we showcase how to employ them to also benefit experimental work.

*The author was partially supported by Vetenskapsrådet Grant 2015-04139 and 2019-03551.

Congruences for Fourier coefficients of modular forms on arithmetic progressions come with two distinct kinds of main applications. There are general Ramanujan-type congruences which are especially important in combinatorics, and the special case of $\mathrm{U}_p$-congruences, which enjoy a much more arithmetic geometric theory via their connection to slopes of $p$-adic modular forms. In the first part of this series, we characterized Ramanujan-type congruences that are not $\mathrm{U}_p$-congruences for integral weight modular forms in terms of congruences satisfied by Hecke eigenvalues. In this work, we extend our investigations to modular forms of half-integral weight. The picture of Ramanujan-type congruences that we arrive at is significantly more subtle than the one for integral weights. Notably, we describe a dichotomy between congruences originating in Hecke eigenvalues and in algebraic parts of twisted central L-values.

A long history of discoveries in combinatorics and number theory originated in Ramanujan's work on divisibility patterns for the partition function [17]:

$$p(5n+4) \equiv 0 \pmod{5}, \quad p(7n+5) \equiv 0 \pmod{7}, \quad p(11n+6) \equiv 0 \pmod{11}.$$

The development set off by Ramanujan bifurcated into a geometric theory of $\mathrm{U}_p$-congruences and a theory of much more combinatorial flavor around Ramanujan-type congruences. Given a weakly holomorphic modular form $f$, a prime $\ell$, and integers $M > 0$ and $\beta$, we say that $f$ has a Ramanujan-type congruences modulo $\ell$ on the arithmetic progression $M\mathbb{Z} + \beta$ if its Fourier coefficients $c(f;\, n)$ satisfy

$$\forall n \in \mathbb{Z} : c(f;\, Mn + \beta) \equiv 0 \pmod{\ell}.$$

We say that this congruence is maximal, if $f$ does not satisfy a Ramanujan-type congruence modulo $\ell$ on $M'\mathbb{Z} + \beta$ for any proper divisor $M'$ of $M$. Ramanujan's original congruences lie at the intersection of $\mathrm{U}_p$-congruences and general Ramanujan-type congruences, since $M$ in his cases is a power of $\ell$.

Given the rich and powerful theory around $\mathrm{U}_p$-operators and geometric modular forms, it is highly desirable to provide a stronger connection between them and Ramanujan-type congruences. This is what we achieve in the present paper. Our first theorem asserts that Ramanujan-type congruences are preserved by the action of the shallow Hecke algebra. Since Atkin-Lehner operators do not preserve all congruences, Theorem A is optimal.

Before we state this theorem, we give four notational remarks. We work with characters of the metaplectic group as opposed to multiplier systems; The character $\chi_\eta$ is the one associated to the Dedekind $\eta$-function, and the subgroups $\widetilde{\Gamma}_0(N) \subseteq \mathrm{Mp}_1(\mathbb{Z})$ of the metaplectic group are the preimages of $\Gamma_0(N)$ under the projection from $\mathrm{Mp}_1(\mathbb{R})$ to $\mathrm{SL}_2(\mathbb{R})$. We write the Fourier expansion of a modular form $f$ associated with $\chi_\eta^r$, $\gcd(r, 24) = 1$, as $\sum c(f; n) \exp(2\pi i\, n\tau/24)$, thus discarding the denominator 24. We write $\mathrm{T}_p^{\mathrm{cl}}$ for the Hecke operator on half-integral weight modular forms that is associated with double-cosets of determinant $p^2$. Finally, for a rational prime $\ell$, $\ell$-integral numbers are the elements of $\mathbb{Q}$ with $\ell$-free denominator.

**Theorem A.** *Fix a prime $\ell$, a weight $k \in \frac{1}{2} + \mathbb{Z}$, and a character $\chi\chi_\eta^r$, $\gcd(r, 24) = 1$, of $\widetilde{\Gamma}_0(N)$, where $\chi$ is a Dirichlet character modulo $N$. Fix a positive integer $M$ and an* *integer $\beta$. Then the space of cusp forms $f$ of weight $k$ for $\chi\chi_\eta^r$ with $\ell$-integral Fourier coefficients $c(f;\,n)$ that satisfy the Ramanujan-type congruence*

$$\forall n \in \mathbb{Z} : c(f;\, Mn+\beta) \equiv 0 \pmod{\ell}$$

*is stable under the action of the Hecke operators $\mathrm{T}_p^{\mathrm{cl}}$, $p \nmid \gcd(M, \ell N)$.*

*Remark.* In Theorem A as in the whole introduction, we restrict our attention to modular forms for the eta-character with rational Fourier coefficients, but in the body of the paper we give statements for modular forms for both the eta-character and the theta-character, and with Fourier coefficients in arbitrary number fields, where $\ell$ then is a prime ideal in the corresponding ring of integers. We refer the reader to the proofs of the main statements in Section 4 for references to the more general results, which require more notation.

In the case of integral weight modular forms, we showed that Ramanujan-type congruences are explained completely by congruences of Hecke eigenvalues. The case of half-integral weights is more intricate, reflecting the appearance of twisted central L-values in the their Fourier coefficients. Our next theorem states a fundamental dichotomy, which highlights the role played by Ramanujan-type congruences on arithmetic progressions $M\mathbb{Z}+\beta$ with cube-free $M$.

We say that a modular form is a generalized Hecke eigenform modulo $\ell$ for $\mathrm{T}_p^{\mathrm{cl}}$ if there is an $\ell$-integral scalar $\lambda_p$ and a nonnegative integer $d$ such that we have the congruence $f|(\mathrm{T}_p^{\mathrm{cl}} - \lambda_p)^{d+1} \equiv 0 \pmod{\ell}$. We refer to $\lambda_p$ as its eigenvalue. Theorem A implies that if $f = \sum_\lambda f_\lambda$ satisfies a Ramanujan-type congruence for generalized Hecke eigenforms $f_\lambda$ for $\mathrm{T}_p^{\mathrm{cl}}$, $p \nmid \gcd(M, \ell N)$, of eigenvalue $\lambda$ that are pairwise different modulo $\ell$, then each $f_\lambda$ satisfies the same Ramanujan-type congruence. In particular, the assumption in Theorem B that $f$ is a generalized Hecke eigenform does not limit its scope, since the canonical idempotents associated with the primary decomposition of the $\mathbb{F}_\ell[\mathrm{T}_p^{\mathrm{cl}}]$-module generated by $f$ yield the required decomposition $\sum_\lambda f_\lambda$ of $f$.

**Theorem B.** *Fix a prime $\ell$, and let $f$ be a cusp form of half-integral weight $k \in \frac{1}{2}+\mathbb{Z}$ for the character $\chi\chi_\eta^r$, $\gcd(r,24)=1$, of $\widetilde{\Gamma}_0(N)$, where $\chi$ is a Dirichlet character modulo $N$, with $\ell$-integral Fourier coefficients $c(f;\,n)$. Assume that $f$ satisfies the maximal Ramanujan-type congruence*

$$\forall n \in \mathbb{Z} : c(f;\, Mn+\beta) \equiv 0 \pmod{\ell}$$

*for some positive integer $M$, $\gcd(M, \ell N) = 1$, and some integer $\beta$. Assume further that $f$ is a generalized Hecke eigenform modulo $\ell$ for all $\mathrm{T}^{\mathrm{cl}}_p$, $p \,|\, M$ prime, of eigenvalue $\lambda_p$.*

*Let $M_{\mathrm{sf}}$ be the largest square-free divisor of $M$, and factor $M = M_{\mathrm{s}}^2 M_{\mathrm{fd}}$ for a positive integer $M_{\mathrm{s}}$ and a square-free, positive integer $M_{\mathrm{fd}}$. Set $\beta_1 = \beta M_{\mathrm{sf}} / M_{\mathrm{s}}^2$. Then we have $\beta_1 \in \mathrm{M}_{\mathrm{fd}} \mathbb{Z}$. Further, if $\ell \ne 2$ set*

$$\widetilde{\mathrm{D}}_k(M\mathbb{Z} + \beta) := \Big\{ D \text{ fund. disc.} : \mathrm{sgn}(D) = (-1)^{k - \frac{1}{2}}, \forall p \,|\, M \text{ prime} : \Big(\frac{|D|}{p}\Big) = \Big(\frac{\beta_1}{p}\Big) \Big\},$$

*and if $\ell = 2$ set*

$$\widetilde{\mathrm{D}}_k(M\mathbb{Z} + \beta) := \Big\{ D \text{ fund. disc.} : \mathrm{sgn}(D) = (-1)^{k - \frac{1}{2}},\ \forall p \,|\, M \text{ prime} \Big(\frac{|D|}{p}\Big) = \Big(\frac{\beta_1}{p}\Big), \\ \forall p \,|\, M_{\mathrm{fd}} : \Big(\frac{|D|/p}{p}\Big) = \Big(\frac{\beta_1/p}{p}\Big) \Big\}.$$

*Then the eigenvalues $\lambda_p$ modulo $\ell$ of $f$ satisfy one of the specific algebraic congruences given in Proposition 3.12, or we have the following Ramanujan-type congruences with cube-free period:*

$$\forall D \in \widetilde{\mathrm{D}}_k(M\mathbb{Z} + \beta) : \forall n \in \mathbb{Z} : c(f;\, M_{\mathrm{sf}} M_{\mathrm{fd}} n + |D|) \equiv 0 \pmod{\ell}.$$

*In the latter case, we in particular have the following congruences on fundamental discriminants:*

$$\forall D \in \widetilde{\mathrm{D}}_k(M\mathbb{Z} + \beta) : c(f;\, |D|) \equiv 0 \pmod{\ell}.$$

*Remark.* (1) As in the case of Theorem A, the analogue of Theorem B holds for algebraic Fourier coefficients.

(2) The scarcity of Ramanujan-type congruences with cube-free period was investigated earlier joint work with Ahlgren and Beckwith [2].

(3) Since the coefficients $c(f;\, |D|)$ of a newform are related to the algebraic parts of twisted central L-values by work of Waldspurger [23] and Kohnen [10], Theorem B describes a dichotomy between Ramanujan-type congruences originating in Hecke eigenvalues and in special L-values. For example, Kohnen's Corollary 1 shows that if the weight-$2k-1$ newform $\tilde{f}$ corresponding to $f$ has level 1, the final congruences in Theorem B are equivalent to

$$\forall D \in \widetilde{\mathrm{D}}_k(M\mathbb{Z}+\beta) : \frac{(k-\frac{3}{2})!\sqrt{|D|}^{2k-2}}{\sqrt{\pi}^{2k-1}} \frac{\langle f, f\rangle}{\langle \tilde{f}, \tilde{f}\rangle} \mathrm{L}(\tilde{f} \otimes \chi_D, 2k-1) \equiv 0 \pmod{\ell},$$

where $\langle\cdot,\cdot\rangle$ is the Petersson scalar product.

(4) If $p \,|\, M_{\mathrm{fd}}$, then the congruences satisfied by generalized Hecke eigenvalues $\lambda_p$ are the same as the ones satisfied in the case of integral weight modular forms. This is not the case if $p \nmid M_{\mathrm{fd}}$. See Proposition 3.12 for details.

As in the case of integral weights, this work relies on a number of structure results that we establish in this paper. It is noteworthy that apart from the square-class result in Statement (1) of the next theorem, which extends previous results by Radu [15, 16], the proofs are completely disjoint from those in the first part of this series.

**Theorem C.** *Fix a prime $\ell$ and let $f$ be a cusp form of half-integral weight $k \in \frac{1}{2}+\mathbb{Z}$ for the character $\chi\chi_\eta^r$, $\gcd(r,24)=1$, of $\widetilde{\Gamma}_0(N)$, where $\chi$ is a Dirichlet character modulo $N$, with $\ell$-integral Fourier coefficients $c(f;n)$. Assume that $f$ satisfies the Ramanujan-type congruence*

$$\forall n \in \mathbb{Z} : c(f; Mn+\beta) \equiv 0 \pmod{\ell}$$

*for some positive integer $M$ and some integer $\beta$. Consider a prime $p \,|\, M$ that does not divide $\ell N$ and factor $M$ as $M_p M_p^\#$ with a $p$-power $M_p$ and $M_p^\#$ co-prime to $p$.*

*(1) We have the congruence on square-classes*

$$\forall n \in \mathbb{Z},\, n \,\square\, \beta \pmod{M} : c(f; n) \equiv 0 \pmod{\ell},$$

*where we write $n \,\square\, \beta \pmod{M}$ if $n$ and $\beta$ lie in the same square-class modulo $M$, that is, if there is an integer $u$, $\gcd(u,M)=1$, with $n \equiv u^2\beta \pmod{M}$.*

*(2) We have the Ramanujan-type congruence*

$$\forall n \in \mathbb{Z} : c(f; M'n+\beta) \equiv 0 \pmod{\ell} \quad \text{with } M' = \gcd(M, M_{\mathrm{sf}}\beta),$$

*where $M_{\mathrm{sf}}$ is the largest square-free divisor of $M$, if $M$ is odd, and in general $M_{\mathrm{sf}}$ is defined as* $\gcd(8,M)(M_2^\#)_{\mathrm{sf}}$.

*(3) If $M_p \,|\, \beta$, we have the Ramanujan-type congruence*

$$\forall n \in \mathbb{Z} : c(f;\, M^{\#}_p n + \beta) \equiv 0 \pmod{\ell}.$$

*(4) Assume that the Ramanujan-type congruence of $f$ modulo $\ell$ on $M\mathbb{Z} + \beta$ is maximal and that $M_p$ is a square. If $M_p = p^2$ assume further that $\ell \ne 2$. Then we have the Ramanujan-type congruence with gap*

$$\forall n \in \mathbb{Z} \setminus p\mathbb{Z} : c(f;\, (M/p)n + M_p \beta') \equiv 0 \pmod{\ell} \quad \text{with } M_p \beta' \equiv \beta \pmod{M^{\#}_p}.$$

*(5) If $M$ is co-prime to $\ell N$ and $f$ is a Hecke eigenform modulo $\ell$ for all Hecke operators $\mathrm{T}^{\mathrm{cl}}_q$, $q \,|\, M$ prime, then there is $p \,|\, M$ such that*

$$\forall n \in \mathbb{Z} : c(f;\, M_p n + \beta) \equiv 0 \pmod{\ell}.$$

*Remark.* For ease of notation, we provide Theorem C only for rational Fourier coefficients. It also holds for prime ideals $\ell$ and eigenforms with Fourier coefficients in the ring of $\ell$-integers of a number field. See the proof of Theorem C in Section 4.

As remarked in the first paper of the present series [18], Ono's and Ahlgren–Ono's congruences for the partition function [3, 14] differ significantly from Atkin's earlier ones [4] in so far as after fixing $\beta \pmod{M_\ell}$ the former naturally arise on arithmetic progressions $M\mathbb{Z} + \beta$ for two square-classes of $\beta \pmod{M_p}$, $p \ne \ell$, while the latter occur for a single square-class of $\beta \pmod{M_p}$. Recent work by Ahlgren–Allen–Tang [1] ensures the existence of infinitely many Ramanujan-type congruences like Atkin's, and Johansson has found many of them computationally [8]. Statement (4) in Theorem C provides an explanation for the two square-classes of $\beta \pmod{M_p}$ in the setting of Ono, because then $M_p$ is a fourth power when inspecting Ramanujan-type congruences of, for instance, $\mathrm{U}_\ell(1/\eta)$. It does not, however, apply to Atkin's congruences, for which $M_p$ is a cube. The next theorem asserts that two square-classes of $\beta \pmod{M_p}$ cannot occur for any Ramanujan-type congruence alike the ones that Atkin discovered (see [24] for further examples of such congruences). In this sense, Statement (4) in Theorem C is sharp.

**Theorem D.** *Fix a prime $\ell$ and let $f$ be a cusp form of half-integral weight $k \in \frac{1}{2} + \mathbb{Z}$ for the character $\chi \chi^r_\eta$, $\gcd(r, 24) = 1$, of $\widetilde{\Gamma}_0(N)$, where $\chi$ is a Dirichlet character modulo $N$, with $\ell$-integral Fourier coefficients $c(f;\, n)$. Assume that $f$ satisfies the maximal Ramanujan-type congruence*

$$\forall n \in \mathbb{Z} : c(f;\, p^2(pn + \beta_0)) \equiv 0 \pmod{\ell}$$

*for some prime $p \nmid \ell N$ and some integer $\beta_0$, $p \nmid \beta_0$.*

*Then for an integer $\beta'_0$, we have the Ramanujan-type congruence*

$$\forall n \in \mathbb{Z} : c(f;\, p^2(pn + \beta'_0)) \equiv 0 \pmod{\ell} \quad \textit{if and only if} \quad \Big(\frac{\beta'_0}{p}\Big) = \Big(\frac{\beta_0}{p}\Big).$$

*Remark.* As for the preceding theorems, there is an analogue of Theorem C for algebraic Fourier coefficients.

We next discuss some consequences for partition congruences. One of the most persistent obstacles in the theory of Ramanujan-type partition congruences are potential congruences of, say, the generating series $\mathrm{U}_\ell(1/\eta)$ to very sparse series, akin to theta series. They can be connected to Ramanujan-type congruences on $M\mathbb{Z} + \beta$ for cube-free $M$, which were the theme of previous joint work with Ahlgren and Beckwith [2]. The next theorem provides a further obstruction of algebraic flavor to such congruences.

To appreciate the essence of this obstruction, it is helpful to recall three facts. First, a theorem by Bruinier [5] on Hecke eigenforms guarantees under rather mild conditions the nonvanishing modulo $\ell$ of Fourier coefficients at fundamental discriminants subject to prescribed square-class conditions at finitely many primes. Second, not every Hecke eigenform modulo $\ell$ is also a Hecke eigenform. Third, not every congruence among newforms of integral weight yields congruences among the corresponding Shintani lifts. In fact in the beginning of the theory it was not clear whether congruences between newforms of half-integral weight exist at all; Maeda [12] answered a corresponding question by Hida in 1983. Only recently, it was demonstrated by Dummigan [6] that such congruences are in fact rather ubiquitous.

**Theorem E.** *Fix a prime $\ell > 3$ and write $p(n)$ for the partition function. Assume that we have the maximal Ramanujan-type congruence*

$$\forall n \in \mathbb{Z} : p\Big(\frac{24\ell M n + \beta}{24}\Big) \equiv 0 \pmod{\ell}, \quad \delta := \Big(\frac{1-\beta}{\ell}\Big),$$

*for some positive, cube-free integer $M \neq 1$, $\ell \nmid M$, and some integer $\beta$.*

*Set $k = (\ell - 2)/2$ if $\delta = 0$, and $k = (\ell^2 - 2)/2$ if $\delta = -1$. Then there is a modular form $f$ of weight $k$ for the character $\chi_\eta^{2k}$ of $\mathrm{Mp}_1(\mathbb{Z})$ with the following properties:*

*(i) All Fourier coefficients of $f$ are $\ell$-integral and $f \not\equiv 0 \pmod{\ell}$.*

*(ii) If* $\delta = 0$*, then for all fundamental discriminants* $D \in \widetilde{\mathrm{D}}_k(M\mathbb{Z} + (\beta - 1)/\ell)$ *as in Theorem B, we have a Ramanujan-type congruence for* $f$ *modulo* $\ell$ *on* $M\mathbb{Z} + |D|$*. If* $\delta = -1$*, then we have a Ramanujan-type congruence for* $f$ *modulo* $\ell$ *on* $\ell M\mathbb{Z} + |D|$ *for all* $D \in \widetilde{\mathrm{D}}_k(\ell M\mathbb{Z} + \beta - 1)$*. In particular, for both* $\delta = 0$ *and* $\delta = -1$*, we have the congruence* $c(f; |D|) \equiv 0 \pmod{\ell}$ *for all respective* $D$*.*

We conclude our results with a theorem that dissects the methods of proof applied by Atkin and Ono in their deductions of Ramanujan-type congruences modulo $\ell$ on $\ell p^3\mathbb{Z} + \beta$ and $\ell p^4\mathbb{Z} + \beta$, $\beta \in \frac{-1}{24} + \mathbb{Z}$, for the partition generating function $\eta^{-1}$. In both cases, Hecke operators act diagonally modulo $\ell$ on a suitable space of cusp forms with $\ell$-integral Fourier coefficients. It is natural to ask whether this method of proof, which is omnipresent when deriving any Ramanujan-type congruence, is in fact based on an equivalence.

In preparation of stating Theorem F, given a prime $\ell > 3$ and $\delta \in \{0, -1\}$, we consider cusp forms $f_{\ell,\delta}$ with $\ell$-integral Fourier coefficients as in (1.7) and (1.8) of [2] satisfying

$$\begin{aligned} f_{\ell,0} &\equiv \sum_{n \in \mathbb{Z}} c(\eta^{-1}; \ell n) \exp\big(2\pi i \tfrac{n}{24} \tau\big) \pmod{\ell} \quad \text{and} \\ f_{\ell,-1} &\equiv \sum_{\substack{n \in \mathbb{Z} \\ (\frac{-n}{\ell}) = -1}} c(\eta^{-1}; n) \exp\big(2\pi i \tfrac{n}{24} \tau\big) \pmod{\ell}. \end{aligned}$$

If $\delta = 0$, the weight of $f_{\ell,\delta}$ equals $(\ell - 2)/2$ and its character is $\chi_\eta^{-\ell}$. If $\delta = -1$ its weight equals $(\ell^2 - 2)/2$ and its character is $\chi_\eta^{-1}$. We write $\mathbb{Z}_\ell$ for the ring of $\ell$-integers, and let

$$F_{\ell,\delta} := \mathbb{Z}_\ell\big[\mathrm{T}^{\mathrm{cl}}_p, p > 3 \text{ prime}\big] f_{\ell,\delta}$$

be the $\ell$-integral Hecke module generated by $f_{\ell,\delta}$.

**Theorem F.** *Fix a prime* $\ell > 3$ *and write* $p(n)$ *for the partition function.*

*(1) Assume that we have the maximal Ramanujan-type congruence*

$$\forall n \in \mathbb{Z} : p\Big(\frac{24\ell q^3 n + \beta}{24}\Big) \equiv 0 \pmod{\ell}, \quad \delta := \Big(\frac{1-\beta}{\ell}\Big),$$

*for some prime* $q > 3$*,* $q \neq \ell$*, and some integer* $\beta$*. Define* $\beta_0$ *by* $\beta_0 - 1 = (\beta - 1)/q^2$*, Then* $\beta_0 \in \mathbb{Z}$*. Assume further that we have*

$$\exists n \in \mathbb{Z} : p\Big(\frac{24\ell q n + \beta_0}{24}\Big) \not\equiv 0 \pmod{\ell}.$$

*Then the generalized eigenvalues of* $\mathrm{T}^{\mathrm{cl}}_q$ *acting modulo* $\ell$ *on* $F_{\ell,\delta}$ *are congruent to* $\pm p^{-1}$ *if* $\delta = 0$ *and* $\pm p^{-2}$ *if* $\delta = -1$.

*(2) Assume that we have the maximal Ramanujan-type congruence*

$$\forall n \in \mathbb{Z} : p\Big(\frac{24\ell q^4 n + \beta}{24}\Big) \equiv 0 \pmod{\ell}, \quad \delta := \Big(\frac{1-\beta}{\ell}\Big),$$

*for some prime* $q > 3$, $q \ne \ell$, *and some integer* $\beta$. *Define* $\beta_0$ *by* $\beta_0 - 1 = (\beta - 1)/q^2$. *Then* $\beta_0 \in \mathbb{Z}$. *Assume further that we have*

$$\exists n \in \mathbb{Z} : p\Big(\frac{24\ell q^2 n + \beta_0}{24}\Big) \not\equiv 0 \pmod{\ell}.$$

*Then* $\mathrm{T}^{\mathrm{cl}}_q$ *acting modulo* $\ell$ *on* $F_{\ell,\delta}$ *is nilpotent.*

In order to gauge the potential of Theorem F, we can ask for which $\ell$ we have

$$F_{\ell,0} = \mathrm{S}_{\frac{\ell-2}{2}}\big(\chi_\eta^{-\ell}; \mathbb{Z}_\ell\big) \quad \text{or} \quad F_{\ell,-1} = \mathrm{S}_{\frac{\ell^2-2}{2}}\big(\chi_\eta^{-1}; \mathbb{Z}_\ell\big).$$

In these cases, Theorem F describes the action of $\mathrm{T}^{\mathrm{cl}}_q$ modulo $\ell$ completely, and informs us that the methods of proof employed by Atkin and Ono are indeed essentially (i.e. up to semi-simplification) equivalent to the Ramanujan-type congruences that they derive. In the spirit of the Maeda Conjecture one might be tempted to believe these equalities, but the appearance of $\ell$-integral coefficients as opposed to rational coefficients yield additional complications as illustrated by the partial congruences among newforms that appear in Theorem E.

The results in this paper are primarily enabled by an analysis of the Ramanujan-type congruences satisfied by Shimura lifts of half-integral weight modular forms. Given a half-integral weight modular form $f$, such congruences are governed by the Fourier coefficients $c(f; |D|n^2)$, $n \in M'\mathbb{Z} + \beta'$. As $n$ varies in $M'\mathbb{Z} + \beta'$, the set of integers $|D|n^2$ is not contained in any arithmetic progression $M\mathbb{Z} + \beta$. To relate Ramanujan-type congruences for half-integral weight modular forms to those of integral weight modular forms, we thus need to derive multiplicative congruences, that is, congruences for the Fourier coefficients $c(f; n)$, where $n$ lies in the same square-class modulo $M$ as $\beta$. Paralleling ideas in earlier work of Radu [15, 16], we achieve this in Section 2, in which we employ the representation theoretic tools developed in the first part of this series to prove that congruences on arithmetic progressions $M\mathbb{Z} + \beta$ incur congruences on square-classes modulo $M$. In a specific setting, treated in Section 2.2, we can also show that more than one square-classes modulo $M$ occurs. Note that the representation theoretic statement for $\mathrm{Mp}_1(\mathbb{Z})$ from which we derive this is not true in general. This differs from the case of integral weight modular forms, for which we were able to show in the first part of this series that congruences always occur on two square-classes modulo $M$.

With these square-class results at hand, in Section 3, we proceed to the study of congruences of Shimura lifts. Section 3 does not employ any representation theory, and follows more classical lines of thought. Generalized Hecke eigenforms modulo $\ell$ are the main obstacle to overcome. In Sections 3.2, 3.3, and 3.7, we take ideas from the first paper of this series one step further to understand the effect of almost all Hecke operators on Ramanujan-type congruences in Section 3.8. We have organized Section 3 in many small sections, which each feature one statement of certain independent relevance. In Section 4 we combine all of them to establish the theorems stated in the introduction.

We wish to particularly point to the extended Example 3.19, which illustrates how to explicitly work with Ramanujan-type congruences on $M\mathbb{Z} + \beta$, $M$ not cube-free, even for very large $M$. We base this on Maeda's example of a congruence between modular forms of half-integral weight. Maeda used a prime ideal $\ell$ that is not principal, and in particular, the Fourier coefficients of the modular form that we consider are not rational. Dummigan's results [6] possibly provide a corresponding example with rational coefficients. In light of the historical precedence set by Maeda's example, our Example 3.19 can be viewed as one of the simplest examples possible, yet it reveals some of the challenges that future calculations with Ramanujan-type congruences, say, from additive combinatorics have to meet. This in particular includes Ramanujan-type congruences that appear only at excessively large $M$ that are divisible by more than one prime. A natural combinatorial explanation for a single of these Ramanujan-type congruences could be an important step in the theory, provided that it is indeed genuinely rooted in combinatorics. To contextualize this requirement, recall that all modular forms can be expressed by eta-quotients [9] and therefore afford some (usually artificial) combinatorial interpretation.

*Acknowledgment* The author wishes to thank Scott Ahlgren, Olivia Beckwith, and Olav Richter for inspiring discussions and helpful comments. The author is also grateful to the referee for their helpful and exceptionally thorough report.

# 1 Preliminaries

We adopt large parts of the setup and notation from the first paper in this series [18]. In particular, we will not revisit the theory of $\ell$-kernels developed in Sections 1 and 2 of that work, which were intentionally written in such a way that they subsume the case of the metaplectic group.

We introduce the following symbol for square-classes of integers $n$ and $m$. Given a positive integer $M$, we define

$$m \,\square\, n \pmod{M} \quad \text{if and only if} \quad \exists u \in \mathbb{Z}, \gcd(u, M) = 1 : m = u^2 n.$$

The ring of integers of a number field $K$ will be written as $\mathcal{O}_K$. Given an ideal $\ell \subseteq \mathcal{O}_K$, we write $\mathcal{O}_{K,\ell}$ for the localization of $\mathcal{O}_K$ at $\ell$ and $\mathrm{F}_{K,\ell}$ for its residue field, if $\ell$ is prime. In case of $K = \mathbb{Q}$, we set $\mathbb{Z}_\ell = \mathcal{O}_{K,\ell}$ and moreover identify ideals $\ell \subseteq \mathbb{Z}$ with positive integers.

We write $e(z)$ for $\exp(2\pi i \tau)$ for $\tau \in \mathbb{C}$. The Poincaré upper half plane will be written $\mathbb{H}$. It carries the usual action via Möbius transformations of $\mathrm{SL}_2(\mathbb{R})$, or more generally $\mathrm{GL}_2^+(\mathbb{R})$.

**1.1 The metaplectic group** Recall the general metaplectic group $\mathrm{GMp}_1^+(\mathbb{R})$ and its realization in [18] including the projection to $\mathrm{GL}_2^+(\mathbb{R})$. We have a section of sets

$$\mathrm{GL}_2^+(\mathbb{R}) \hookrightarrow \mathrm{GMp}_1^+(\mathbb{R}),\, \gamma \longmapsto (\gamma, \omega_\gamma), \quad \text{where } \omega_\gamma(i) \in \mathbb{H}. \tag{1.1}$$

We occasionally identify elements of the left hand side with their image under this map.

In order to describe multiplication more explicitly, we write

$$(\gamma, \omega_\gamma)(\delta, \omega_\delta) \;=\; (\gamma\delta, \sigma(\gamma,\delta)\omega_{\gamma\delta}), \tag{1.2}$$

where $\sigma : \mathrm{GMp}_1^+(\mathbb{R}) \times \mathrm{GMp}_1^+(\mathbb{R}) \longrightarrow \{\pm 1\}$ is a 2-cocycle. Its values can be found in work of Maaß [11] or, more explicitly, Strömberg [21]. For the time being, we write $c(\gamma)$ and $d(\gamma)$ for the bottom left and right entries of $\gamma \in \mathrm{GL}_2^+(\mathbb{R})$, and let $\tilde{c}(\gamma) = c(\gamma)$, if $c(\gamma) \ne 0$, and $\tilde{c}(\gamma) = d(\gamma)$, otherwise. Then we have

$$\begin{aligned}
\sigma(\gamma,\delta) &= \big(\tilde{c}(\gamma\delta)\tilde{c}(\gamma),\, \tilde{c}(\gamma\delta)\tilde{c}(\delta)\big)_\infty, && \text{if } c(\gamma), c(\delta), c(\gamma\delta) \ne 0;\\
\sigma(\gamma,\delta) &= \big(\tilde{c}(\gamma),\, \tilde{c}(\delta)\big)_\infty, && \text{if } c(\gamma), c(\delta) \ne 0, c(\gamma\delta) = 0;\\
\sigma(\gamma,\delta) &= \big(\tilde{c}(\gamma),\, \tilde{c}(\delta)\big)_\infty, && \text{if } c(\gamma) = c(\delta) = 0;\\
\sigma(\gamma,\delta) &= \big(\tilde{c}(\gamma),\, -\tilde{c}(\delta)\big)_\infty, && \text{if } c(\gamma) = 0, c(\delta) \ne 0;\\
\sigma(\gamma,\delta) &= \big(-\tilde{c}(\gamma),\, \tilde{c}(\delta)\big)_\infty, && \text{if } c(\gamma) \ne 0, c(\delta) = 0.
\end{aligned}$$

Recall the notation $\widetilde{\Gamma}_0(N)$ for subgroups of $\mathrm{Mp}_1(\mathbb{Z})$. We have two prominent characters of $\widetilde{\Gamma}_0(4)$ and of $\mathrm{Mp}_1(\mathbb{Z})$. Their classical description is based on associated multiply systems

$$\nu_\theta\big(\big(\begin{smallmatrix} a & b \\ c & d \end{smallmatrix}\big)\big) = \big(\tfrac{c}{d}\big)\epsilon_d^{-1},$$

where $\epsilon_d = 1$ if $d \equiv 1 \pmod{4}$ and $\epsilon_d = i$ if $d \equiv 3 \pmod{4}$, and

$$\begin{aligned}
\nu_\eta\big(\big(\begin{smallmatrix} a & b \\ c & d \end{smallmatrix}\big)\big) &= \big(\tfrac{d}{c}\big)e\big(\tfrac{1}{24}\big((a+d)c - bd(c^2-1) - 3c\big)\big), && \text{if } c > 0 \text{ is odd;}\\
\nu_\eta\big(\big(\begin{smallmatrix} a & b \\ c & d \end{smallmatrix}\big)\big) &= \big(\tfrac{c}{d}\big)e\big(\tfrac{1}{24}\big((a+d)c - bd(c^2-1) + 3d - 3 - 3cd\big)\big), && \text{if } c > 0 \text{ is even.}
\end{aligned}$$

Then we define

$$\begin{aligned}
&\chi_\theta : \widetilde{\Gamma}_0(4) \longrightarrow \mathrm{GL}_1(\mathbb{C}), \quad \chi_\theta(1,-1) = -1, \chi_\theta(\gamma, \omega_\gamma) = \nu_\theta(\gamma);\\
&\chi_\eta : \mathrm{Mp}_1(\mathbb{Z}) \longrightarrow \mathrm{GL}_1(\mathbb{C}), \chi_\eta(1,-1) = -1, \chi_\eta(\gamma, \omega_\gamma) = \nu_\eta(\gamma).
\end{aligned}$$

We will also use the notation $\chi_{\theta+}$, which we explain in the next section.

**1.2 Modular forms** Given a weight $k \in \frac{1}{2} + \mathbb{Z}$ and a Dirichlet character $\chi$ modulo $N$, we write $\mathrm{M}_k(\chi\chi_\theta^r)$, $r$ odd, and $\mathrm{M}_k(\chi\chi_\eta^r)$, $\gcd(r, 24) = 1$, for the space of weight-$k$ modular forms on $\widetilde{\Gamma}_0(\mathrm{lcm}(4, N))$ and $\widetilde{\Gamma}_0(N)$ that transform as the given character. When saying that $f$ is a modular form for $\chi\chi_\theta^r$ or $\chi\chi_\eta^r$, where $\chi$ is a Dirichlet character modulo $N$, we always assume that the level of $f$ divides $\mathrm{lcm}(4, N)$ in the former case and $N$ in the latter one. When saying that $f$ is a modular form for $\chi\chi_\theta^r$ or $\chi\chi_\eta^r$, where $\chi$ is a Dirichlet character modulo $N$, we will always assume that the level of $f$ divides $N$. We write $\mathrm{M}_k(N, \chi\chi_\theta^r)$ and $\mathrm{M}_k(N, \chi\chi_\theta^r)$ to emphasize the level. Observe that nonzero modular forms occur only if we have $r \equiv 2k + 1 - \chi(-1) \pmod{4}$.

The subspaces of cusp forms of $\mathrm{M}_k(\chi\chi_\theta^r)$ and $\mathrm{M}_k(\chi\chi_\eta^r)$ will be written as $\mathrm{S}_k(\chi\chi_\theta^r)$ and $\mathrm{S}_k(\chi\chi_\eta^r)$. Even though weakly holomorphic modular forms will make appearance in Section 2, we will not need special notation for them.

In the case of the theta character, we write the Fourier expansion as

$$f(\tau) = \sum_{n \in \mathbb{Z}} c(f; n) e(n\tau) \in \mathrm{M}_k(\chi\chi_\theta^r), \tag{1.3}$$

while in the case of the eta-character, we simplify notation and write

$$f(\tau) = \sum_{n \in \mathbb{Z}} c(f; n) e\big(\tfrac{n}{24}\tau\big) \in \mathrm{M}_k(\chi\chi_\eta^r). \tag{1.4}$$

By slight abuse of notation, we write $\mathrm{M}_k(\chi\chi_{\theta+}^r)$ and $\mathrm{S}_k(\chi\chi_{\theta+}^r)$ for the Kohnen-plus space, even though the associated condition cannot be expressed in terms of a character alone.

Given a number field $K \subset \mathbb{C}$, the spaces of modular forms with Fourier coefficients in $\mathcal{O}_{K,\ell}$ will be written, for instance, as $\mathrm{M}_k(\chi\chi_\theta^r; \mathcal{O}_{K,\ell})$.

Given a arbitrary prime $p$ in the case of the theta-character and $p \nmid 6$ in the case of the eta-character, there are classical Hecke operators, which we denote by $\mathrm{T}_p^{\mathrm{cl}}$ to distinguish them from the vector-valued ones in the first part of this series [18]. In this paper, we can assume that $p$ is odd.

Shimura in Theorem 1.7 of [19] gives an expression for the action of the Hecke operator on a modular form $f \in \mathrm{M}_k(\chi\chi_\theta^r)$ in the case of $r \equiv 2k \pmod{4}$. Note that if $r \equiv 2k+2 \pmod{4}$, i.e., if $\chi(-1) = -1$, we can replace $\chi(p)$ by $\chi(p)\chi_\theta(p)^2 = \chi(p)\chi_4(p)$ to accommodate this assumption. We hence have the formula

$$c\big(f \big| \mathrm{T}_p^{\mathrm{cl}};\, n\big) \;=\; c(f;\, np^2) + \Big(\frac{(-1)^{\frac{r-1}{2}}\, n}{p}\Big)\chi(p)\, p^{k-\frac{3}{2}}\, c(f;\, n) + \chi(p^2)\, p^{2k-2} c(f;\, n/p^2).$$

For a modular form $f \in \mathrm{M}_k(\chi\chi_\eta^r)$, if $p \nmid 6$, we obtain a similar relation:

$$c\big(f \big| \mathrm{T}_p^{\mathrm{cl}};\, n\big) \;=\; c(f;\, np^2) + \Big(\frac{(-1)^{\frac{r-1}{2}}\, 12n}{p}\Big)\chi(p)\, p^{k-\frac{3}{2}}\, c(f;\, n) + \chi(p^2)\, p^{2k-2} c(f;\, n/p^2).$$

Observe that these Hecke operators would in some other papers be referred to as $\mathrm{T}_{p^2}^{\mathrm{cl}}$, since they are associated to double cosets with determinant $p^2$.

**1.3 The Shimura lift** For an integer $D$, which will typical be a fundamental discriminant, we set $\chi_D := (\frac{D}{\bullet})$. We have a map

$$\begin{aligned} \mathrm{V}_{24} :\qquad \mathrm{M}_k\big(\chi v_\eta^r\big) &\longrightarrow \mathrm{M}_k\Big(576N,\, \chi_{12}\chi\chi_\theta^r\Big), \\ \sum_{n \in \mathbb{Z} + \frac{r}{24}} c(n)\, e(n\tau) &\longmapsto \sum_{n \in 24\mathbb{Z} + r} c\big(\tfrac{n}{24}\big)\, e(n\tau). \end{aligned} \tag{1.5}$$

Further, for $k \ge \frac{5}{2}$ we have a family of Shimura lifts [13, 19], which for each fundamental discriminant $D$ with $\mathrm{sgn}(D) = (-1)^{k-\frac{1}{2}}$, yields a map

$$\mathrm{Sh}_D : \mathrm{S}_k\big(\chi\chi_{\theta+}^r\big) \longrightarrow \mathrm{S}_{2k-1}\big(2N, \chi^2\big),$$

which is compatible with Hecke operators $\mathrm{T}_p^{\mathrm{cl}}$, $p \nmid 2N$, and with the property that

$$\sum_{n=1}^{\infty} c\big(f;\, |D| n^2\big) n^{-s} \;=\; \frac{\mathrm{L}\big(\mathrm{Sh}_D(f),\, s\big)}{\mathrm{L}\big(\chi\chi_4^{2k-r}\chi_D,\, s + \frac{3}{2} - k\big)}. \tag{1.6}$$

Note that this is the linearized Shimura lift, which maps newforms to scalar multiples of newforms. Note also that similar to our discussion of Hecke operator following Shimura's work, we extended (1.6) from the case $r \equiv 2k \pmod{4}$ to general $r$.

Combining the congruence $E_{\ell_{\mathbb{Z}}-1} \equiv 1 \pmod{\ell}$ for the Eisenstein series with the Shimura lift we obtain

$$\mathrm{Sh}_D : \mathrm{S}_{k \,(\mathrm{mod}\, \ell_{\mathbb{Z}}-1)}\big(\chi \chi_{\theta+}^r\big) \pmod{\ell} \longrightarrow \mathrm{S}_{2k-1 \,(\mathrm{mod}\, \ell_{\mathbb{Z}}-1)}\big(2N, \chi^2\big) \pmod{\ell}, \tag{1.7}$$

for which (1.6) holds coefficient-wise modulo $\ell$ and where $\ell_{\mathbb{Z}}\mathbb{Z} = \ell \cap \mathbb{Z}$ is the prime ideal above $\ell$.

To obtain a Shimura lift for modular forms of the eta-character, we compose $\mathrm{Sh}_D$ with $\mathrm{V}_{24}$. We keep the same notation, since the meaning is unambiguous thanks to the requirement that $\gcd(r, 24) = 1$. We have a maps

$$\begin{aligned} \mathrm{Sh}_D : \mathrm{S}_k\big(\chi \chi_\eta^r\big) &\longrightarrow \mathrm{S}_{2k-1}\big(288N, \chi^2\big),\, f \longmapsto \mathrm{Sh}_D(\mathrm{V}_{24}\, f), \\ \mathrm{Sh}_D : \mathrm{S}_{k \,(\mathrm{mod}\, \ell_{\mathbb{Z}}-1)}\big(\chi \chi_\eta^r\big) \pmod{\ell} &\longrightarrow \mathrm{S}_{2k-1 \,(\mathrm{mod}\, \ell_{\mathbb{Z}}-1)}\big(288N, \chi^2\big) \pmod{\ell}. \end{aligned} \tag{1.8}$$

For the sake of clarity, we note that the target space of (1.8) can be narrowed down as in the work of Yang on the case of the trivial character and level one [25]. Ahlgren communicated to us that a forthcoming joint work with Andersen and Dicks will treat this generalization.

To contextualize our results, we occasionally appeal to the Shintani lift [20]. It yields modular forms of half-integral weight starting from modular forms of integral weight. Up to a proportionality factor it is an inverse to the Shimura lifts.

**1.4 Congruences** We quote Theorem 2.2, Proposition 2.4, and Corollary 2.5 of [18] for clarity and later reference. We refer the reader to [18] for notation, and remark that we will apply the statement for $(V, \phi) = (\mathbb{C}\, f, \mathrm{id})$, where $f$ is a weakly holomorphic modular form.

**Theorem 1.1.** *Let $K \subset \mathbb{C}$ be a number field with fixed complex embedding, and let $\ell$ be an ideal in $\mathcal{O}_K$. Let $(V, \phi)$ be an abstract module of weakly holomorphic modular forms for $\Gamma \subseteq \mathrm{Mp}_1(\mathbb{Z})$ realized in weight $k \in \frac{1}{2}\mathbb{Z}$. Assume that* $\ker_\Gamma(V)$ *is a congruence subgroup of level $N$. Let $N_\ell$ be the smallest positive integer such that*

$$\gcd\big(\ell, N / N_\ell\big) = 1 \text{ if } k \in \mathbb{Z} \quad \text{and} \quad \gcd\big(\ell, \mathrm{lcm}(4, N) / N_\ell\big) = 1 \text{ if } k \in \tfrac{1}{2} + \mathbb{Z}. \tag{1.9}$$

*Then the $\ell$-kernel of $(V, \phi)$ is a right-module for $\mathcal{O}_{K,\ell}[\widetilde{\Gamma}_0(N_\ell) \cap \Gamma]$.*

**Proposition 1.2.** *Fix a positive integer $M$, a rational number $\beta$, and an abstract module of weakly holomorphic modular forms $(V, \phi)$ for $\Gamma \subseteq \mathrm{Mp}_1(\mathbb{Z})$, $T \in \Gamma$, over a ring that contains $e(\beta/M)$. Fix a $T$-eigenvector $v \in V$ with eigenvalue $e(\beta)$. Then we have*

$$(\mathrm{T}_M \phi)\big(\mathrm{T}_M(v, \beta)\big) \;=\; M^{1-\frac{k}{2}} \sum_{n \in \beta + M\mathbb{Z}} c(f;\, n) e\big(\tfrac{n}{M}\tau\big). \tag{1.10}$$

**Corollary 1.3.** *Fix a positive integer $M$, a rational number $\beta$, a number field $K \subset \mathbb{C}$ that contains $e(\beta/M)$, and an ideal $\ell \subseteq \mathscr{O}_K$ that is co-prime to $M$. Consider an abstract module of weakly holomorphic modular forms $(V, \phi)$ for $\Gamma \subseteq \mathrm{Mp}_1(\mathbb{Z})$, $T \in \Gamma$, over $\mathscr{O}_{K,\ell}$. Fix a $T$-eigenvector $v \in V$ with eigenvalue $e(\beta)$. Then the Ramanujan-type congruence*

$$\forall n \in \mathbb{Z} : c(f;\, Mn + \beta) \equiv 0 \pmod{\ell}$$

*for the modular form $f = \phi(v)$ is equivalent to*

$$\mathrm{T}_M(v, \beta) \in \ker_{\mathrm{FE}\ell}\big(\mathrm{T}_M(V, \phi)\big).$$

# 2 Multiplicative congruences

As opposed to the case of integral weights, our investigation of Ramanujan-type congruences in half-integral weights does not require much input from modular representation theory, and could in fact be formulated without. We gather all of it in this section. We start the section by clarifying the role of square-classes modulo $M$ in relation to Ramanujan-type congruences on arithmetic progressions $M\mathbb{Z} + \beta$. This parallels ideas of Radu [15, 16], but the mechanics of our proof are greatly streamlined by appealing to the representation theoretic perspective. Section 2.2 depends on that perspective in a more essential way. If $M$ is exactly divisible by the square of a prime, we establish a coupling result for two square-classes modulo $M$. This yields Ramanujan-type congruences with gaps, which were dominant in previous work [18], but do not always occur in the half-integral weight setting. Finally in Section 2.3, we provide an essential tool to control the support of Fourier coefficients modulo $\ell$.

**2.1 Square-classes** There is an intimate connection between the action of $\Gamma_0(N)$ and $\widetilde{\Gamma}_0(N)$ and square-classes modulo $M$ in the arithmetic progressions $M\mathbb{Z} + \beta$. The next proposition is a corner stone of our discussion in Section 3.

For clarity, note that $\chi$ in the statement is not necessarily associated with a Dirichlet character, but can be any character of $\Gamma_0(N)$. Since we assume in the proposition that the kernel of $\chi$ is a congruence subgroup of level $N_\chi$, a weakly holomorphic modular form $f$ for $\chi\chi_\theta^r$ or $\chi\chi_\eta^r$ has a Fourier expansion of the shape

$$\sum_{n\in\mathbb{Z}} c\big(f;\tfrac{n}{N_\chi}\big)\, e\big(\tfrac{n}{N_\chi}\tau\big) \quad \text{or} \quad \sum_{n\in\mathbb{Z}} c\big(f;\tfrac{n}{N_\chi}\big)\, e\big(\tfrac{n}{\mathrm{lcm}(24,N_\chi)}\tau\big).$$

**Proposition 2.1.** *Fix a number field $K\subset\mathbb{C}$ and a prime ideal $\ell\subset\mathcal{O}_K$. Let $f$ be a weakly holomorphic modular form of half-integral weight for a character $\chi\chi_\theta^r$, $r$ odd, or $\chi\chi_\eta^r$,* $\gcd(r,24)=1$, *with Fourier coefficients $c(f;\,n)\in\mathcal{O}_{K,\ell}$, where $\chi$ is character of $\Gamma_0(N)$ for some positive integer $N$ whose kernel is a congruence subgroup. Assume that $f$ satisfies the Ramanujan-type congruence*

$$\forall n\in\mathbb{Z} : c(f;\,Mn+\beta)\equiv 0 \pmod{\ell}$$

*for some positive integer $M$ and some $\beta\in N_\chi^{-1}\mathbb{Z}$, where $N_\chi$ is the level of the kernel of $\chi\chi_\theta^r$ or $\chi\chi_\eta^r$, respectively.*

*Then for all $\beta'=u^2\beta$, $u\in\mathbb{Z}$, with* $\gcd(u,MN_\chi)=1$, *we have the Ramanujan-type congruence*

$$\forall n\in\mathbb{Z} : c(f;\,Mn+\beta')\equiv 0 \pmod{\ell}.$$

*In particular, if $\beta$ is integral, then we have*

$$\forall n\in\mathbb{Z} : c(f;\,M'n+\beta)\equiv 0 \pmod{\ell},$$

*where $M'=\gcd(M,M_{\mathrm{sf}}\beta)$ and $M_{\mathrm{sf}}=\gcd(8,M)\prod_p p$ with $p$ running through all odd prime divisors of $M$. Further, we then have*

$$\forall n\in\mathbb{Z},\, n\,\square\,\beta \pmod{M'} : c(f;\,n)\equiv 0 \pmod{\ell}.$$

*Proof.* We can and will assume that $K$ contains all $M$-th roots of unity and $M^{\frac{k}{2}}$. The second part follows from the first one in the same way as in the proof of Proposition 3.5 of [18]. To establish the first part, we consider the vector $M^{\frac{k}{2}-1}\mathrm{T}_M(f,\beta)$ in the abstract module of modular forms $\mathrm{T}_M\mathbb{C}f$. We have to show that $M^{\frac{k}{2}-1}\mathrm{T}_M(f,\beta')$ lies in $\ker_{\mathrm{FE}\ell}(\mathrm{T}_M\mathbb{C}f)$ to then apply Proposition 1.2 and obtain the result. As all operations in this proof are defined over the ring of cyclotomic integers, we can discard the normalizing factor $M^{\frac{k}{2}-1}$.

Observe that we have $N\,|\,N_\chi$. Since only $\beta' \pmod{M}$ is relevant in the statement, we can replace $N$ by any multiple of it and modify $u$, if needed, by adding a suitable multiple of $M$ to ensure that it is co-prime to $N_\chi$. In particular, we can and will assume that $24 \,|\, N$. Let $\overline{u} \in \mathbb{Z}$ denote the inverse of $u$ modulo $MN_\chi$. We fix a matrix $\gamma \in \Gamma_0(N)$ that is congruent to $\left(\begin{smallmatrix} u & 0 \\ 0 & \overline{u} \end{smallmatrix}\right)$ modulo $MN_\chi$. We want to show that applying $(\gamma, \omega_\gamma)$ to $\mathrm{T}_M(f, \beta)$ yields a multiple of $\mathrm{T}_M(f, \beta')$. We can employ the same argument as in the integral weight case, except that we have to track the precise action by scalars of the center of $\mathrm{Mp}_1(\mathbb{Z})$ and the contribution of $\chi_\theta$ or $\chi_\eta$.

For any $h \in \mathbb{Z}$, we find in perfect analogy with the case of integral weight that

$$\delta := \begin{pmatrix} 1 & h \\ 0 & M \end{pmatrix} \gamma \begin{pmatrix} 1 & \overline{u}^2 h \\ 0 & M \end{pmatrix}^{-1} \in \Gamma_0(N_\chi).$$

Using (1.1), we identify each factor in the definition of $\delta$ with the corresponding element in $\mathrm{GMp}_1(\mathbb{Z})$. Then the explicit description of the cocycle in (1.2) shows that the product yields $(\delta, \omega_\delta)$, since the bottom left entries of the first and last factor vanish and their bottom right entries are positive.

In the case of the theta-character, we have

$$\begin{aligned} \mathrm{T}_M(f,\beta)\big|(\gamma,\omega_\gamma) &= \sum_{h \pmod{M}} f \otimes e\big(-\tfrac{h\beta}{M}\big)\, \mathfrak{u}_{\left(\begin{smallmatrix} 1 & h \\ 0 & M \end{smallmatrix}\right)} \big|(\gamma,\omega_\gamma) \\ &= \sum_{h \pmod{M}} \chi\chi_\theta(\delta,\omega_\delta)\, f \otimes e\big(-\tfrac{h\beta}{M}\big)\, \mathfrak{u}_{\left(\begin{smallmatrix} 1 & \overline{u}^2 h \\ 0 & M \end{smallmatrix}\right)}. \end{aligned}$$

Observe that we suppress the dependency of $\delta$ on $h$ from our notation. We finish the proof if we show that $\chi\chi_\theta((\delta,\omega_\delta))$ is independent of $h$. Similarly, in the case of the eta-character, we have to show that $\chi\chi_\eta((\delta,\omega_\delta))$ is independent of $h$.

We have $\delta \equiv \left(\begin{smallmatrix} u & 0 \\ 0 & \overline{u} \end{smallmatrix}\right) \pmod{N_\chi}$, which shows that the image of $\delta$ under the level-$N_\chi$ character $\chi$ is independent of $h$.

It remains to inspect the images under $\chi_\theta$ and $\chi_\eta$, whose kernels satisfy

$$\begin{aligned} \ker \chi_\theta &= \big\{ \big(\gamma, \nu_\theta(\gamma)\omega_\gamma\big) : \gamma \in \Gamma_1(4) \big\}, \\ \ker \chi_\eta &\supseteq \big\{ \big(\gamma, \nu_\eta(\gamma)\omega_\gamma\big) : \gamma \in \Gamma_0^0(12) \big\}, \end{aligned}$$

where $\Gamma_0^0(12) \subseteq \Gamma_0(12)$ consists of matrices whose top right entry is divisible by 12. Recall that we are assuming that $24 \,|\, N$. Using quadratic reciprocity, we find that the image of $(\delta,\omega_\delta)$ under $\chi_\theta$ is constant with respect to $h$, since the bottom left entry of $\delta$ is independent of $h$ and even. To show that $\chi_\eta((\delta,\omega_\delta))$ is independent of $h$, we use the explicit formula for the eta-character, and merely need to assert that $e(d/8)$ is constant with respect to $h$, where $d$ is the bottom left entry of $\delta$. This is the case, since $d \equiv \overline{u} \pmod{8}$. ∎

The proof of Proposition 2.1 can be adapted to establish the next result, which allows us to focus on subgroups $\widetilde{\Gamma}_0(N)$ as opposed to $\widetilde{\Gamma}_1(N)$.

**Lemma 2.2.** *Fix a Dirichlet character $\chi$ of level $N$, and integers $M > 0$ and $\beta$. Let $f$ be a weakly holomorphic modular form of half-integral weight $k$ for a character $\chi\chi_\theta^r$, $r$ odd. Then the Fourier expansion*

$$\sum_{\substack{n \in \mathbb{Z} \\ n \square \beta \,(\mathrm{mod}\, M)}} c(f;\, n)\, e(n\tau)$$

*is associated with a weakly holomorphic modular form of weight $k$ for the character $\chi\chi_\theta^r$ of $\widetilde{\Gamma}_0(M^2 \mathrm{lcm}(N,4))$.*

*Let $f$ be a weakly holomorphic modular form of half-integral weight $k$ for a character $\chi\chi_\eta^r$, $\gcd(r, 24) = 1$. Then the Fourier expansion*

$$\sum_{\substack{n \in \mathbb{Z} \\ n \square \beta \,(\mathrm{mod}\, M)}} c(f;\, n)\, e\big(\tfrac{n}{24}\tau\big)$$

*is associated with a weakly holomorphic modular form of weight $k$ for the character $\chi\chi_\eta^r$ of $\widetilde{\Gamma}_0(M^2 N)$.*

*Proof.* We can and will assume that $K$ contains all $M$-th roots of unity. Consider the case of the theta-character. We can express the given Fourier expansion as

$$\sum_{\substack{\beta' \square \beta \,(\mathrm{mod}\, M) \\ h \,(\mathrm{mod}\, M)}} e\big(\tfrac{-h\beta'}{M}\big) f\big|_k \big(\begin{smallmatrix} M & h \\ 0 & M \end{smallmatrix}\big).$$

Given any matrix $\gamma = \big(\begin{smallmatrix} a & b \\ c & d \end{smallmatrix}\big) \in \Gamma_0(M^2 \mathrm{lcm}(4, N))$, we set

$$\delta \;=\; \big(\begin{smallmatrix} M & h \\ 0 & M \end{smallmatrix}\big)\, \gamma \, \big(\begin{smallmatrix} M & hd\overline{a} \\ 0 & M \end{smallmatrix}\big) \;=\; \begin{pmatrix} a + \frac{hc}{M} & b + \frac{hd - ha\overline{a}d + h^2 \overline{a}dc/M}{M} \\ c & d - \frac{h\overline{a}dc}{M} \end{pmatrix},$$

where $\overline{a}$ is an inverse of $a$ modulo $M^2$. Since $4 \,|\, c$, we can use the same argument as in the proof of Proposition 2.1 to show that $f|\delta = f|\gamma$, confirming the first case of the lemma.

To handle the second case, we can use the expression

$$\sum_{\substack{\beta' \square \beta \,(\mathrm{mod}\, M) \\ h \,(\mathrm{mod}\, M)}} e\big(\tfrac{-h\beta'}{24M}\big) f\big|_k \big(\begin{smallmatrix} M & h \\ 0 & M \end{smallmatrix}\big),$$

where an extra denominator 24 in the exponential twist arises from our rescaling of the Fourier indices in (1.4). Proposition 2.1 allows us to assume that $\gcd(M,6) = 1$. Then, since $h$ in the previous sum is well-defined modulo $M$, we can and will assume that $24 \,|\, h$.

The expression for $\delta$ in terms of the entries of $\gamma \in \Gamma_0(M^2N)$ remains the same as in the case of the theta-character. If we have $\gamma = \left(\begin{smallmatrix} 1 & \pm 1 \\ 0 & 1 \end{smallmatrix}\right)$ or $\gamma = \left(\begin{smallmatrix} -1 & 0 \\ 0 & -1 \end{smallmatrix}\right)$, we directly verify that $\chi_\eta(\delta) = \chi_\eta(\gamma)$. For $\gamma$ with $c > 0$, we can evaluate $\chi_\eta(\delta)$ using the formulas in Section 1.1. Any $\gamma$ with $c < 0$ can be written as a product of matrices which we have already checked. ∎

**2.2 Ramanujan-type congruences with gap** In the case of integral weight we have seen in previous work [18] that square-classes of congruences occur in pairs. This is not the case for all Ramanujan-type congruences in half-integral weights. The next statement shows that in some cases this coupling of square-classes can nevertheless be encountered.

**Proposition 2.3.** *Fix a number field $K \subset \mathbb{C}$ and a prime ideal $\ell \subset \mathcal{O}_K$, $\ell \nmid 2$. Let $f$ be a weakly holomorphic modular form of half-integral weight for a character $\chi\chi_\theta^r$, $r$ odd, with Fourier coefficients $c(f;n) \in \mathcal{O}_{K,\ell}$, where $\chi$ is a Dirichlet character modulo $N$. Assume that $f$ satisfies the Ramanujan-type congruence*

$$\forall n \in \mathbb{Z} : c(f; Mn + \beta) \equiv 0 \pmod{\ell}$$

*for some positive integer $M$ and some integer $\beta$.*

*Given a prime $p$, we factor $M$ as $M_p M_p^\#$ with a $p$-power $M_p$ and $M_p^\#$ co-prime to $p$. For some $p \,|\, M$, $\gcd(p, \ell N)$, assume that $M_p = p^2$ and that $p$ exactly divides $\beta$. Then for all integers $\beta'$ such that $p$ exactly divides $\beta'$ and $\beta' \equiv \beta \pmod{M_p^\#}$, we have*

$$\forall n \in \mathbb{Z} : c(f; Mn + \beta') \equiv 0 \pmod{\ell}.$$

The proof of Proposition 2.3 requires an indivisibility result for Kloosterman sums

$$\mathrm{K}_2(\psi, a) := \sum_{x_1 x_2 \equiv a \pmod{p}} \psi(x_1 + x_2).$$

Here, $\psi$ is a multiplicative character of period $p$ for a prime $p$, and $a$ is an integers. The author is grateful to Will Sawin, who pointed out a proof of the following lemma.

**Lemma 2.4.** *Let $\ell \neq 2$ and $p$ be distinct primes, $a$ an integer, and $\psi$ a nontrivial multiplicative character of period $p$. Then we have*

$$\mathrm{K}_2(\psi, a) \not\equiv 0 \pmod{\ell}.$$

*Remark 2.5.* If $\ell = 2$ and $a$ is not a square modulo $p$, the conclusion of the lemma does not hold. This matches the situation in the proof of Proposition 2.3, which hence requires the assumption $\ell \nmid 2$.

*Proof of Lemma 2.4.* Recall the equality $\mathrm{K}_2(\psi^b, a) = \mathrm{K}_2(\psi, ab^2)$ for nonzero integers $b$ modulo $p$. It allows us to assume that $\psi(x) = e(x/p)$.

The Kloosterman sum takes values in algebraic integers of the $p$-th cyclotomic field:

$$\mathbb{Z}(\zeta_p) \cong \mathbb{Z}[X] / \big( X^{p-1} + X^{p-2} + \cdots + 1 \big) \cong \big( \mathbb{Z}[X] / (X^p - 1) \big) \big/ \big( X^{p-1} + X^{p-2} + \cdots + 1 \big).$$

Reduction modulo $\ell$ intertwines with the polynomial quotient, yielding a ring over the finite field $\mathbb{F}_\ell$. It therefore suffices to examine the sum

$$\sum_{\substack{x \,(\mathrm{mod}\, p) \\ x\overline{x} \equiv 1 \,(\mathrm{mod}\, p)}} X^{ax + \overline{x}} = \sum_{i=0}^{p-1} c_i X^i \in \mathbb{Z}[X] / (X^p - 1).$$

The lemma follows if we show that the coefficients $c_i$ of $X^i$, $0 \le i < p$, are not all equal modulo $\ell$. We suppose the opposite and derive a contradiction.

If $a \equiv 0 \pmod{p}$, then $c_i = 1$ for $0 < i < p$, but $c_0 = 0$, finishing the proof. We can therefore assume that $a$ is nonzero modulo $p$.

We notice that $ax + \overline{x} \equiv h \pmod{p}$ for any $h \pmod{p}$ yields a quadratic equation in $x$, which has at most two solutions. We thus have $c_i \in \{0, 1, 2\}$. Since $\ell$ is odd, 0, 1, and 2 are distinct modulo $\ell$. We are assuming that all $c_i$ are equal modulo $\ell$, and hence they coincide as integers. Comparing the number of terms in the sum, with the possible values of the $c_i$, we find $p - 1 = \sum_i c_i = p c_0 \in \{0, p, 2p\}$, a contradiction. ■

*Proof of Proposition 2.3.* If $\ell \mid 3$, then every prime satisfies $p \equiv \pm 1 \pmod{\ell}$ and the statement is vacuous. We can therefore assume that $\ell \nmid 6$. The statement is also trivial if $p = 2$, since there is only one nonzero square-class modulo 2. To simplify notation, we will assume that $4 \mid N$.

We can and will assume that $K$ contains all $M$-th roots of unity. Proposition 2.1 and Lemma 2.2 allow us to replace $f$ by

$$\sum_{\substack{n \in \mathbb{Z} \\ n \,\square\, \beta \,(\mathrm{mod}\, M_p^\#)}} c(f;\, n) e(\tau n).$$

In other words, we may assume that $M = p^2$. Using Corollary 1.3, our assumptions yield $\mathrm{T}_M(f,\beta) \in \ker_{\mathrm{FE}\ell}(\mathrm{T}_M \mathbb{C} f)$, and we need to show that $\mathrm{T}_M(f,\beta')$ occurs in that $\ell$-kernel for all integers $\beta'$ that are exactly divisible by $p$. By Proposition 2.1 this is the case for all $\beta'$ that lie in the same square-class modulo $M$ as $\beta$. Therefore, we can and will assume that $\beta'$ lies in the opposite square-class of $\beta$ modulo $M$.

We start by inspecting the eigenspaces of $(T,\omega_T) \in \widetilde{\Gamma}_0(N)$ acting on $\ker_{\mathrm{FE}\ell}(\mathrm{T}_M \mathbb{C} f)$. Notice that in the next computation, all matrices are upper triangular, which allows us to suppress the cocycles in (1.2) from our considerations. Using the usual set of representatives, we have for $0 \le b < M$

$$f \otimes \mathfrak{u}_{\left(\begin{smallmatrix} M & 0 \\ 0 & 1 \end{smallmatrix}\right)} \big| \left(\begin{smallmatrix} 1 & 1 \\ 0 & 1 \end{smallmatrix}\right) = f \big| \left(\begin{smallmatrix} 1 & M \\ 0 & 1 \end{smallmatrix}\right) \otimes \mathfrak{u}_{\left(\begin{smallmatrix} M & 0 \\ 0 & 1 \end{smallmatrix}\right)}, \qquad f \otimes \mathfrak{u}_{\left(\begin{smallmatrix} p & 0 \\ 0 & p \end{smallmatrix}\right)} \big| \left(\begin{smallmatrix} 1 & 1 \\ 0 & 1 \end{smallmatrix}\right) = f \big| \left(\begin{smallmatrix} 1 & 1 \\ 0 & 1 \end{smallmatrix}\right) \otimes \mathfrak{u}_{\left(\begin{smallmatrix} p & 0 \\ 0 & p \end{smallmatrix}\right)},$$
$$f \otimes \mathfrak{u}_{\left(\begin{smallmatrix} 1 & b \\ 0 & M \end{smallmatrix}\right)} \big| \left(\begin{smallmatrix} 1 & 1 \\ 0 & 1 \end{smallmatrix}\right) = f \otimes \mathfrak{u}_{\left(\begin{smallmatrix} 1 & b+1 \\ 0 & M \end{smallmatrix}\right)}, \qquad f \otimes \mathfrak{u}_{\left(\begin{smallmatrix} 1 & M-1 \\ 0 & M \end{smallmatrix}\right)} \big| \left(\begin{smallmatrix} 1 & 1 \\ 0 & 1 \end{smallmatrix}\right) = f \big| \left(\begin{smallmatrix} 1 & 1 \\ 0 & 1 \end{smallmatrix}\right) \otimes \mathfrak{u}_{\left(\begin{smallmatrix} 1 & 0 \\ 0 & M \end{smallmatrix}\right)}.$$

We conclude that the eigenspaces with $T$-eigenvalues $e(\beta/M)$ and $e(\beta'/M)$ are one-dimensional, since $\beta$ and $\beta'$ are not divisible by $M$. These eigenspaces are spanned by $\mathrm{T}_M(f,\beta)$ and $\mathrm{T}_M(f,\beta')$.

We next examine the action of a transformation $(\gamma,\omega_\gamma) \in \widetilde{\Gamma}_0(N)$ that is congruent to $\left(\begin{smallmatrix} 1 & 0 \\ 1 & 1 \end{smallmatrix}\right)$ modulo $M$ on $\mathrm{T}_M(f,\beta) \in \ker_{\mathrm{FE}\ell}(\mathrm{T}_M \mathbb{C} f)$. For an integer $h$ with $h \not\equiv -1 \pmod{p}$, we have

$$\left(\begin{smallmatrix} 1 & h \\ 0 & M \end{smallmatrix}\right) \gamma \;\equiv\; \left(\begin{smallmatrix} 1+h & h \\ M & M \end{smallmatrix}\right) \;\equiv\; \delta \begin{pmatrix} 1 & h\overline{(1+h)} \\ 0 & M \end{pmatrix} \pmod{M},$$

where $\overline{1+h}$ is an inverse modulo $M$ of $1+h$ and $\delta \in \Gamma_0(N)$. We have

$$(\delta,\omega_\delta) \;=\; \left(\begin{smallmatrix} 1 & h \\ 0 & M \end{smallmatrix}\right) (\gamma,\omega_\gamma) \begin{pmatrix} 1 & h\overline{(1+h)} \\ 0 & M \end{pmatrix}^{-1} \in \widetilde{\Gamma}_0(N).$$

As in the proof of Proposition 2.1, we suppress the dependency of $\delta$ on $h$ from our notation.

We use the explicit formula for the cocycle in (1.2) to calculate that

$$f \otimes \mathfrak{u}_{\left(\begin{smallmatrix} 1 & h \\ 0 & M \end{smallmatrix}\right)} \Big| (\gamma,\omega_\gamma) \;=\; f \Big| (\delta,\omega_\delta) \otimes \mathfrak{u}_{\left(\begin{smallmatrix} 1 & h\overline{(1+h)} \\ 0 & M \end{smallmatrix}\right)}.$$

When inserting the definition of $\mathrm{T}_M(f,\beta)$, we obtain

$$\mathrm{T}_M(f,\beta) \big| (\gamma,\omega_\gamma) \;=\; \sum_{h \,(\mathrm{mod}\, M)} e\Big(\frac{-\beta h}{M}\Big) f \otimes \mathfrak{u}_{\left(\begin{smallmatrix} 1 & h \\ 0 & M \end{smallmatrix}\right)} \Big| (\gamma,\omega_\gamma)$$
$$= \sum_{\substack{h \,(\mathrm{mod}\, M) \\ h \not\equiv -1 \,(\mathrm{mod}\, p)}} e\Big(\frac{-\beta h}{M}\Big) f | (\delta,\omega_\delta) \otimes \mathfrak{u}_{\left(\begin{smallmatrix} 1 & h\overline{(1+h)} \\ 0 & M \end{smallmatrix}\right)} + \sum_{\substack{h \,(\mathrm{mod}\, M) \\ h \equiv -1 \,(\mathrm{mod}\, p)}} e\Big(\frac{-\beta h}{M}\Big) f \otimes \mathfrak{u}_{\left(\begin{smallmatrix} 1 & h \\ 0 & M \end{smallmatrix}\right)} \Big| (\gamma,\omega_\gamma). \tag{2.1}$$

We record that all terms in the second summand are associated with matrices of the form $\left(\begin{smallmatrix} p & * \\ 0 & p \end{smallmatrix}\right)$ or $\left(\begin{smallmatrix} M & * \\ 0 & 1 \end{smallmatrix}\right)$. Their contributions later in the proof will vanish.

We next inspect the action of $(\delta, \omega_\delta)$ on $f$, which transforms like $\chi \chi_\theta^r$, where $\chi$ is a Dirichlet character. Let $c$ and $d$ be the bottom entries of $\gamma$. Then the bottom entries of $\delta$ equal $p^2 c$ and $d - h\overline{(1+h)}c$. Recall that the value of $\chi_\theta(\gamma)$ is given by a Jacobi symbol and by $\epsilon_d$. Since we have $N \,|\, c$ and $\chi$ is a Dirichlet character modulo $N$, we conclude that $\chi(d - h\overline{(1+h)}c) = \chi(d)$. By quadratic reciprocity we have

$$\Big(\frac{p^2 c}{d - h\overline{(1+h)}c}\Big) = \Big(\frac{c/c'^2}{d - h\overline{(1+h)}c}\Big) = \Big(\frac{c/c'^2}{d}\Big) = \Big(\frac{c}{d}\Big)$$

for the maximal square divisor $c'^2$ of $c$. Since $4 \,|\, N$ and $N \,|\, c$, the $\epsilon$-values in the formula for $\chi_\theta((\gamma, \omega_\gamma))$ and $\chi_\theta((\delta, \omega_\delta))$ coincide by $d - h\overline{(1+h)}c \equiv d \pmod{4}$. We conclude that the action of $(\delta, \omega_\delta)$ on $f$ is independent of $h$. It only depends on $\gamma$. Hence we can define a complex scalar $\tilde{\omega}_\gamma$ by

$$f \big| (\delta, \omega_\delta) =: \tilde{\omega}_\gamma f.$$

We next apply the following projection to the expression in (2.1):

$$\mathrm{T}_M(f, \beta)\big|_k (\gamma, \omega_\gamma)\Big| \sum_{h' \pmod{M}} e\Big(\frac{-\beta' h'}{M}\Big) \left(\begin{smallmatrix} 1 & h' \\ 0 & 1 \end{smallmatrix}\right).$$

A calculation reveals that the image has eigenvalue $e(\beta'/M)$ under $T$. We have already determined that the corresponding eigenspace is spanned by $\mathrm{T}_M(f, \beta')$. Since the action of $\left(\begin{smallmatrix} 1 & h' \\ 0 & 1 \end{smallmatrix}\right)$ preserves matrices of the form $\left(\begin{smallmatrix} p & * \\ 0 & p \end{smallmatrix}\right)$ and $\left(\begin{smallmatrix} M & * \\ 0 & 1 \end{smallmatrix}\right)$, the terms in the second sum in (2.1) with $h \equiv -1 \pmod{p}$ vanish under this projection. We use the fact that $\chi\chi_\theta(T) = 1$, to find that the remaining terms yield

$$\begin{aligned} \tilde{\omega}_\gamma \sum_{h' \pmod{M}} \sum_{\substack{h \pmod{M} \\ h \not\equiv -1 \pmod{p}}} e\Big(\frac{-\beta h - \beta' h'}{M}\Big) f \otimes \mathfrak{u}_{\left(\begin{smallmatrix} 1 & \frac{h\overline{(1+h)}+h'}{M} \\ 0 & \end{smallmatrix}\right)} \\ = \tilde{\omega}_\gamma \mathrm{T}_M(f, \beta') \sum_{\substack{h, h' \pmod{M} \\ h \not\equiv -1 \pmod{p} \\ h\overline{(1+h)}+h' \equiv 0 \pmod{M}}} e\Big(\frac{-\beta h - \beta' h'}{M}\Big). \end{aligned}$$

The last factor on the right hand side is a multiple of a Kloosterman sum, which we can identify more explicitly by replacing $h$ by $h-1$ and realizing that the condition $h\overline{(1+h)} + h' \equiv 0 \pmod{M}$ becomes $\overline{h} - 1 \equiv h' \pmod{M}$, where $\overline{h}$ is an inverse of $h$ modulo $M$. This leads us to

$$\sum_{\substack{h,h' \pmod{M} \\ h \not\equiv -1 \pmod{p} \\ h\overline{(1+h)}+h' \equiv 0 \pmod{M}}} e\Big(\frac{-\beta h - \beta' h'}{M}\Big) = e\Big(\frac{\beta}{M}\Big) \sum_{\substack{h \pmod{M} \\ h \not\equiv 0 \pmod{p}}} e\Big(\frac{-\beta h - \beta'(\overline{h}-1)}{M}\Big)$$

$$= e\Big(\frac{\beta+\beta'}{M}\Big) \sum_{h \pmod{M}^\times} e\Big(\frac{(-\beta/p)h + (-\beta'/p)\overline{h}}{p}\Big) = e\Big(\frac{\beta+\beta'}{M}\Big) p \,\mathrm{K}_2(\psi_p^b, a),$$

where $\psi_p(h) = e(h/p)$, $a = \beta\beta'/p^2$, and $b = -\beta'/p$. Recall the equality $\mathrm{K}_2(\psi_p^b, a) = \mathrm{K}_2(\psi_p, ab^2)$.

To finish the proof, we have to show that for $a = \beta\beta'/p^2$, and $b = -\beta'/p$ as before the Kloosterman sum $\mathrm{K}_2(\psi_p, ab^2)$ does not vanish modulo $\ell$. Proposition 2.1 allows us to replace $\beta$ by $\beta h^2$ for any integer $h$ that is not divisible by $p$. It therefore suffices to show that

$$\exists h \in \mathbb{Z} \setminus p\mathbb{Z} : \mathrm{K}_2\big(\psi_p,\, h^2 \beta\beta'^3/p^4\big) \not\equiv 0 \pmod{\ell}.$$

By contraposition, we assume that for all such $h$ the congruence holds.

In the remainder of the proof, we can replace $K$ by the $p$-th cyclotomic field, in which the Kloosterman sum takes its value. We identify the Galois group of $K$ with the units in $\mathbb{Z}/p\mathbb{Z}$. Under this identification, the action of $h \in \mathbb{Z} \setminus p\mathbb{Z} \longrightarrow \mathbb{Z}/p\mathbb{Z}$ is given by $\mathrm{K}_2(\psi, a) \longmapsto \mathrm{K}_2(\psi, h^2 a)$. We combine this with

$$\forall h \in \mathbb{Z} \setminus p\mathbb{Z} : \mathrm{K}_2\big(\psi,\, h^2 \beta\beta'^3/p^4\big) \equiv 0 \pmod{\ell}$$

to conclude that

$$\mathrm{K}_2\big(\psi,\, \beta\beta'^3/p^4\big) \equiv 0 \pmod{\mathrm{Nm}(\ell)},$$

where $\mathrm{Nm}(\ell) \in \mathbb{Z}$ is the norm of $\ell$. By Lemma 2.4 this is a contradiction, finishing our proof. ∎

**2.3 The support of Fourier expansions modulo $\ell$** The final statement in this section allows us to control oldforms modulo $\ell$. In the case of integral weights, we could provide a much improved version thanks to the kindness of Serre, who communicated an argument directly leveraging the theory of Galois representations associated with modular forms. We reproduce this statement, as it will be employed several times.

**Lemma 2.6 (Serre).** *Fix a number field $K \subset \mathbb{C}$ and a prime ideal $\ell \subset \mathscr{O}_K$. Let $f$ be a modular form of integral weight $k$ for a Dirichlet character $\chi$ modulo $N$ with Fourier coefficients $c(f;\,n) \in \mathscr{O}_{K,\ell}$. Assume that there is an integer $n \ne 0$ with the property that $c(f;\,n) \not\equiv 0 \pmod{\ell}$. Then given a density-zero set $P$ of primes that are co-prime to $\ell N$, there is $n \in \mathbb{Z}$ with $p \nmid n$ for all $p \in P$ and $c(f;\,n) \not\equiv 0 \pmod{\ell}$.*

Since the connection between Fourier coefficients and Hecke eigenvalues is much weaker in the half-integral weight case than the integral weight case, the reasoning to obtain Lemma 2.6 does not extend, and we fall back to an argument from modular representation theory.

**Lemma 2.7.** *Fix a number field $K \subset \mathbb{C}$ and a prime ideal $\ell \subset \mathscr{O}_K$. Let $f$ be a weakly holomorphic modular form of half-integral weight for a character $\chi\chi_\theta^r$, $r$ odd, or $\chi\chi_\eta^r$,* $\gcd(r,24) = 1$, *with Fourier coefficients $c(f;\,n) \in \mathscr{O}_{K,\ell}$, where $\chi$ is a Dirichlet character modulo $N$.*

*Assume that $f$ satisfies the congruences*

$$\forall n \in \mathbb{Z}, \gcd(n,M) = 1 : c(f;\,n) \equiv 0 \pmod{\ell},$$

*for some positive integer $M$ that is co-prime to $\ell 4N$. Then $f$ is congruent to a constant modulo $\ell$.*

*In particular, consider a linearly independent set of cusp forms $f_i \pmod{\ell}$ of half-integral weight for a character $\chi\chi_\theta^r$, $r$ odd, or $\chi\chi_\eta^r$,* $\gcd(r,24) = 1$, *with Fourier coefficients $c(f_i;\,n) \in \mathscr{O}_{K,\ell}$. Then there is a finite set of integers $n$ co-prime to $M$ such that the matrix with entries $c(f_i;\,n)$ is invertible modulo $\ell$.*

*Proof.* The second part follows from the first one by mere linear algebra. We focus on the first statement. In the case of the $\eta$-character, we have $c(f;\,n) = 0$ if $\gcd(n,6) \ne 1$, which allows us to assume that $\gcd(M,6) = 1$.

For both the eta- and the theta-character, given a prime $p \mid M$ and $\beta \in \mathbb{Z}$ with $\beta$ co-prime to $M_p^\#$, we can combine Proposition 2.1 and Lemma 2.2 to replace $f$ by the modular form

$$f_p^\#(\tau) = \sum_{\substack{n \in \mathbb{Z} \\ n \,\square\, \beta \pmod{M_p^\#}}} c(f;n)e(n\tau).$$

We have

$$\forall n \in \mathbb{Z}, \gcd(n,p) = 1 : c(f_p^\#;\,n) \equiv 0 \pmod{\ell}.$$

Assuming that the lemma is true in the case of $M$ prime, we can conclude that $f_p^\#$ is congruent to a constant modulo $\ell$. In other words, we have

$$\forall n \in \mathbb{Z}, \gcd(n, M_p^\#) = 1 : c(f;\, n) \equiv 0 \pmod{\ell}.$$

Using induction on the number of prime divisors of $M$, we obtain the lemma, provided that we can establish the case of $M$ prime. In the remainder of the proof, we assume $M = p$ is prime.

We consider the $\ell$-kernel of $\mathrm{T}_p \mathbb{C} f$. By our assumptions and Corollary 1.3, we find that

$$\forall \beta \in \mathbb{Z} \setminus p\mathbb{Z} : \mathrm{T}_p(f, \beta) \in \ker_{\mathrm{FE}\ell}\big(\mathrm{T}_p \mathbb{C} f\big).$$

We fix a matrix $\gamma \in \Gamma_0(\mathrm{lcm}(4, N))$ that is congruent to $\left(\begin{smallmatrix} 1 & 0 \\ 1 & 1 \end{smallmatrix}\right)$ modulo $p$. As in Equation (2.1) in the proof of Proposition 2.3, we find that $\mathrm{T}_p(f,\beta)|(\gamma,\omega_\gamma)$ equals

$$\sum_{\substack{h \pmod{p} \\ h \not\equiv -1 \pmod{p}}} e\Big(\frac{-\beta h}{p}\Big) f|(\delta, \omega_\delta) \otimes \mathfrak{u}_{\left(\begin{smallmatrix} 1 & h\overline{(1+h)} \\ 0 & p \end{smallmatrix}\right)} + e\Big(\frac{\beta}{p}\Big) f\big|(\delta', \omega_{\delta'}) \otimes \mathfrak{u}_{\left(\begin{smallmatrix} p & 0 \\ 0 & 1 \end{smallmatrix}\right)}, \tag{2.2}$$

where $\overline{1+h}$ is an inverse modulo $p$ of $1+h$ and

$$\delta = \left(\begin{smallmatrix} 1 & h \\ 0 & p \end{smallmatrix}\right) (\gamma, \omega_\gamma) \begin{pmatrix} 1 & h\overline{(1+h)} \\ 0 & p \end{pmatrix}^{-1} \quad \text{and} \quad \delta' = \left(\begin{smallmatrix} 1 & h \\ 0 & p \end{smallmatrix}\right) (\gamma, \omega_\gamma) \left(\begin{smallmatrix} p & 0 \\ 0 & 1 \end{smallmatrix}\right)^{-1}.$$

Similar to the projection applied in the proof of Proposition 2.3, since $\ell \nmid p$, we can project onto the subspace that is invariant under $(T, \omega_T) \in \mathrm{Mp}_1(\mathbb{Z})$, discarding all contributions from the vectors $\mathrm{T}_p(f,\beta)$, $\beta \not\equiv 0 \pmod{p}$. Since the last term in (2.2) is invariant under $(T, \omega_T)$, we obtain a non-zero linear combination

$$f \otimes \Big(c_0 \mathrm{T}_p(f, 0) + c_1 \mathfrak{u}_{\left(\begin{smallmatrix} p & 0 \\ 0 & 1 \end{smallmatrix}\right)}\Big) \in \ker_{\mathrm{FE}\ell}\big(\mathrm{T}_p \mathbb{C} f\big).$$

The coefficients $c_0$ and $c_1$ are $\ell$-integral and $c_1$ is invertible modulo $\ell$.

We next examine the associated Fourier expansion. After extending $K$ by $\sqrt{p}$, if needed, we obtain that

$$c_0 \sum_{n \in p\mathbb{Z}} c(f;n) e\Big(\frac{n}{p}\tau\Big) + c_1 p^{\frac{k}{2}} \sum_{n \in \mathbb{Z}} c(f;n) e(pn\tau) \equiv 0 \pmod{\ell}$$

The result is immediate if $c_0 \equiv 0 \pmod{\ell}$. Otherwise, we proceed as in the proof of Theorem 3.1 (3.5) in [18], by induction on the power of $p$ exactly dividing $n$. ∎

# 3 Ramanujan-type and Hecke congruences

We fix a setup, which we will refer to repeatedly throughout this Section 3. Fix a number field $K \subset \mathbb{C}$ and a prime ideal $\ell \subset \mathcal{O}_K$. Consider a cusp form

$$f \in \mathrm{S}_k\big(\chi\chi_{\theta+}^r; \mathcal{O}_{K,\ell}\big) \quad \text{or} \quad f \in \mathrm{S}_k\big(\chi\chi_{\eta}^r; \mathcal{O}_{K,\ell}\big) \tag{3.1}$$

of half-integral weight $k$ for one of the characters $\chi\chi_{\theta+}^r$, $r$ odd, or $\chi\chi_{\eta}^r$, $\gcd(r,24) = 1$, with Fourier coefficients $c(f;\, n) \in \mathcal{O}_{K,\ell}$, where $\chi$ is a Dirichlet character modulo $N$. To simplify notation, we assume that $4 \,|\, N$ in the case of the theta-character. For clarity, we emphasize that we work with the Kohnen-plus space associated in Section 1.2 with $\chi_{\theta+}$ as opposed to $\chi_\theta$.

We consider a positive integer $M$, $\gcd(M, \ell N) = 1$, and an integer $\beta$ for which we have a Ramanujan-type congruence

$$\forall n \in \mathbb{Z} :\; c(f;\, Mn + \beta) \equiv 0 \pmod{\ell}. \tag{3.2}$$

We will often assume that this Ramanujan-type congruence is maximal.

Proposition 2.1 says that $f$ has a Ramanujan-type congruence on the arithmetic progression $\gcd(M, M_{\mathrm{sf}}\beta)\mathbb{Z} + \beta$. If the Ramanujan-type congruence in (3.2) is maximal, we therefore have $M \,|\, M_{\mathrm{sf}}\beta$.

Since $M$ and $N$ are co-prime, $M$ is odd in the case of that $f$ is a cusp form for $\chi\chi_{\theta+}^r$. If $f$ is a cusp form for $\chi\chi_{\eta}^r$, then we have $c(f;n) \equiv 0 \pmod{\ell}$ if $n \not\equiv r \pmod{24}$. Therefore if $\gcd(M, \beta, 6) \ne 1$, then the Ramanujan-type congruence on $M\mathbb{Z} + \beta$ is trivial. In other words, every maximal and nontrivial Ramanujan-type congruence (3.2) satisfies $\gcd(M, 6) = 1$. We assume throughout that every Ramanujan-type is nontrivial in this sense.

Given a prime $p$, we factor $M$ as $M_p M_p^\#$ with a $p$-power $M_p$ and $M_p^\#$ co-prime to $p$. We set $M_{\mathrm{sf}} = \gcd(8, M) \prod_p p$ with $p$ running through all odd prime divisors of $M$.

**3.1 Multiplicative congruences** In order to conveniently apply the results about congruences on square-classes modulo $M$, we factorize

$$M\mathbb{Z} + \beta \;=\; M_1 (M_0 \mathbb{Z} + \beta_0), \quad \text{where} \quad M_1 = \gcd(M, \beta),\, M_0 = M / \gcd(M, \beta),\, \beta_0 = \beta / \gcd(M, \beta). \tag{3.3}$$

A slight refinement of this is given by

$$M\mathbb{Z} + \beta \;=\; M_{\mathrm{fd}} M_{\mathrm{s}}^2 (M_0 \mathbb{Z} + \beta_0), \quad \text{where } M_{\mathrm{fd}} M_{\mathrm{s}}^2 = \gcd(M, \beta), \tag{3.4}$$

for a positive integer $M_{\mathrm{s}}$ and square-free $M_{\mathrm{fd}}$.

With this notation and by Proposition 2.1, the Ramanujan-type congruence (3.2) implies the multiplicative congruence

$$\forall n_0 \in \mathbb{Z}, n_0 \square \beta_0 \pmod{M_0} : c(f; M_1 n_0) \equiv 0 \pmod{\ell}. \tag{3.5}$$

**3.2 Decomposition into generalized Hecke eigenforms** In this section, we fix notation for the decomposition of modular forms modulo $\ell$ of half-integral weight into generalized Hecke eigenforms.

As in the case of integral weights, generalized Hecke eigenforms arise from congruences among newforms. Since the Shimura lift preserves integrality of Fourier coefficients, but the Shintani lift does not, it can happen that modular forms in half-integral weights are not congruent modulo $\ell$, while all their Shimura lifts are proportional modulo $\ell$.

**Example 3.1.** Using, for instance the computer algebra system Sage [22], we arrive at the following. We have a congruence modulo 3 between the level 11 and level 33 newforms

$$f_1(\tau) = e(\tau) - 2e(2\tau) - e(3\tau) + 2e(4\tau) + \mathcal{O}\big(e(5\tau)\big),$$
$$f_2(\tau) = e(\tau) + e(2\tau) - e(3\tau) - e(4\tau) + \mathcal{O}\big(e(5\tau)\big).$$

This yields a generalized Hecke eigenform $(f_1 - f_2)/3$ modulo 3. The Shintani lifts of $f_1$ and $f_2$ yield modular forms of weight $\frac{3}{2}$ that are proportional to

$$g_1(\tau) = e(3\tau) + e(4\tau) + \mathcal{O}\big(e(11\tau)\big),$$
$$g_2(\tau) = e(3\tau) + \mathcal{O}\big(e(11\tau)\big).$$

Now $g_1$ and $g_2$ are not congruent modulo 3, but all their Shimura lifts are integral multiplies of $f_1$ and $f_2$, which are congruent to each other.

Consider a modular form $f$ as in (3.2). After extending $K$ if necessary, for every prime $p \,|\, M$, the $\mathrm{F}_{K,\ell}[\mathrm{T}^{\mathrm{cl}}_p]$-module generated by $f$ is a direct sum of primary modules $\mathrm{F}_{K,\ell}[\mathrm{T}^{\mathrm{cl}}_p] / (\mathrm{T}^{\mathrm{cl}}_p - \lambda_p)^{d_p+1}$ by the primary decomposition. Applying to $f$ the polynomials in $\mathrm{T}^{\mathrm{cl}}_p$ that are the associated canonical idempotents, we can decompose $f$ as a sum $\sum_\lambda f_\lambda$ of generalized Hecke eigenforms $f_\lambda$ modulo $\ell$ of eigenvalues $\lambda = (\lambda_p)_{p \,|\, M}$ for the classical Hecke operators $\mathrm{T}^{\mathrm{cl}}_p$, $p \,|\, M$. For $t = (t_p)_{p \,|\, M}$ with $0 \le t_p \le d_p$, we set

$$f_{\lambda,t} := f_\lambda \Big| \prod_{p \,|\, M} (\mathrm{T}^{\mathrm{cl}}_p - \lambda_p)^{t_p}. \tag{3.6}$$

Given $\lambda$ as before, we write $d_p = d_{\lambda,p}$ for some nonnegative integer with the property that $f_\lambda | (\mathrm{T}_p^{\mathrm{cl}} - \lambda_p)^{d_p+1} \equiv 0 \pmod{\ell}$. We define the span of all $f_{\lambda,t}$ as

$$F_\lambda \;:=\; \mathrm{span}\, \mathrm{F}_{K,\ell}\big\{ f_{\lambda,t} : 0 \le t_p \le d_p \text{ for all } p \,|\, M \big\}.$$

Its $\mathrm{F}_{K,\ell}$-linear dual $F_\lambda^\vee$ and the dual of its Shimura lifts $\mathrm{Sh}_D(F_\lambda)^\vee$ will appear later.

In case of the theta-character, the $D$-th Shimura lift is Hecke-equivariant away from the level $N$, and in case of eta-character away from $6N$. In particular, for every fundamental discriminant $D$ with $\mathrm{sgn}(D) = (-1)^{k-\frac{1}{2}}$ and every prime $p \,|\, M$, we have

$$\mathrm{Sh}_D(f_\lambda) \big| \big(\mathrm{T}_p^{\mathrm{cl}\mathbb{Z}} - \lambda_p\big)^{d_p+1} \equiv 0 \pmod{\ell}, \tag{3.7}$$

where $\mathrm{T}_p^{\mathrm{cl}\mathbb{Z}}$ is the classical Hecke operator in integral weight, associated with matrices of determinant $p$ as opposed to $p^2$. Observe that $\mathrm{Sh}_D(f_\lambda) \pmod{\ell}$ might be annihilated by a lower power of $\mathrm{T}_p^{\mathrm{cl}\mathbb{Z}} - \lambda_p$ or even vanish.

**3.3 Decomposition of L-series for Shimura lifts** As a next step, we provide a decomposition of L-series associated with Shimura lifts. In particular, we show how Ramanujan-type congruences as in (3.2) are reflected by congruences of individual terms in this decomposition.

We start by recalling Lemma 4.10 of [18]. The condition in that Lemma that $d$, which corresponds to our $d_\lambda$ in Section 3.2, is minimal is for convenience only, avoiding redundancies. We will account for such redundancies in a different way and may hence discard that assumption. The meaning of the subscript of $f_t$ differs also, in order to better fit our discussion of several as opposed to a single prime $p$.

In the next proposition, the L-series away from an integer $M$ and the L-polynomial at a prime $p$ of a modular form $\tilde{f}$ of integral weight $k$ and Dirichlet character $\chi$ will be written as

$$\mathrm{L}^M(\tilde{f}, s) \;:=\; \sum_{\substack{n=1 \\ \gcd(n,M)=1}}^{\infty} c(\tilde{f}; n) n^{-s}, \quad \mathrm{L}_p\big(\lambda_p, X\big) \;:=\; 1 - \lambda_p X + \chi(p) p^{k-1} X^2.$$

Observe that we suppress from our notation the dependency of the L-polynomial on $k$ and $\chi$. We define congruences between (formal) Dirichlet series as congruences between their coefficients.

In the next lemma as in (3.7), the Hecke operators $\mathrm{T}_p^{\mathrm{cl}\mathbb{Z}}$ are the classical Hecke operators in integral weight, and hence associated to matrices of determinant $p$, not $p^2$.

**Lemma 3.2 (Lemma 4.10 of [18]).** *Fix a number field $K \subset \mathbb{C}$ and a prime ideal $\ell \subset \mathscr{O}_K$. Let $\tilde{f}$ be a modular form of integral weight $k$ for a Dirichlet character $\chi$ modulo $N$ with Fourier coefficients $c(\tilde{f};\, n) \in \mathscr{O}_{K,\ell}$. Consider a prime $p$ that is co-prime to $\ell N$. Assume that we have the congruence*

$$\tilde{f} \big| \big(\mathrm{T}_p^{\mathrm{cl}\mathbb{Z}} - \lambda_p\big)^{d+1} \equiv 0 \pmod{\ell}, \quad \lambda_p \in \mathscr{O}_{K,\ell},$$

*for some nonnegative integer $d$. For any nonnegative integer $t$, set*

$$\tilde{f}_t := \tilde{f} \big| \big(\mathrm{T}_p^{\mathrm{cl}\mathbb{Z}} - \lambda_p\big)^{t}.$$

*Then we have*

$$\mathrm{L}(\tilde{f}, s) \equiv \sum_{t=0}^{d} \frac{p^{-st}\, \mathrm{L}^{p}(\tilde{f}_t, s)}{\mathrm{L}_p(\lambda_p, p^{-s})^{t+1}} \pmod{\ell}. \tag{3.8}$$

To state the next corollary, recall the notation $f_{\lambda,t}$ set up in Section 3.2.

**Corollary 3.3.** *Fix a level-$N$ cusp form $f$ and a prime ideal $\ell$ as in (3.1), and a positive integer $M$ with* $\gcd(M, \ell N) = 1$. *Then we have*

$$\mathrm{L}\big(\mathrm{Sh}_D(f), s\big) \equiv \sum_{\lambda} \sum_{t=0}^{d_\lambda} \mathrm{L}^{M}\big(\mathrm{Sh}_D(f_{\lambda,t}), s\big) \prod_{p \mid M} \frac{p^{-st_p}}{\mathrm{L}_p\big(\lambda_p, p^{-s}\big)^{t_p+1}} \pmod{\ell}. \tag{3.9}$$

*Remark 3.4.* Explicitly, the L-polynomial in this corollary equals

$$\mathrm{L}_p\big(\lambda_p, X\big) = 1 - \lambda_p X + \chi^2(p) p^{2k-2} X^2,$$

since the weight and character of $\mathrm{Sh}_D(f)$ are $2k-1$ and $\chi^2$.

*Proof.* For convenience, we may assume that $M$ is square-free. We prove the statement by induction on the number of prime divisors of $M$. The case that $M = p$ is prime follows directly from Lemma 3.2.

Given any $M$ and a prime $q \nmid M$, observe that Lemma 3.2 also implies

$$\mathrm{L}^{M}\big(\mathrm{Sh}_D(f), s\big) \equiv \sum_{\lambda_q} \sum_{t_q=0}^{d_{\lambda,q}} \mathrm{L}^{qM}\big(\mathrm{Sh}_D(f_{\lambda_q,t_q}), s\big) \frac{q^{-st_q}}{\mathrm{L}_q\big(\lambda_q, q^{-s}\big)^{t_q+1}} \pmod{\ell}.$$

Keeping the shorthand notation $\lambda = (\lambda_p)_{p \mid M}$ and $t = (t_p)_{p \mid M}$, we will need the relation

$$\big(f_{\lambda,t}\big)_{\lambda_q,t_q} \equiv f_{\lambda',t'} \pmod{\ell} \quad \text{with } \lambda' = (\lambda_p)_{p \mid qM},\ t' = (t_p)_{p \mid qM},$$

which follows directly from the definition.

Now assuming that the lemma holds for given $M$, we can insert the above relation into (3.9) and obtain the desired congruence from

$$\mathrm{L}\big(\mathrm{Sh}_D(f), s\big) \equiv \sum_{\lambda} \sum_{t}^{d_\lambda} \mathrm{L}^M\big(\mathrm{Sh}_D(f_{\lambda,t}),\, s\big) \prod_{p \mid M} \frac{p^{-s t_p}}{\mathrm{L}_p\big(\lambda_p, p^{-s}\big)^{t_p+1}}$$
$$\equiv \sum_{\lambda} \sum_{t}^{d_\lambda} \sum_{\lambda_q} \sum_{t_q=0}^{d_{\lambda,q}} \mathrm{L}^{qM}\big(\mathrm{Sh}_D(f_{\lambda',t'}),\, s\big) \frac{q^{-s t_q}}{\mathrm{L}_q\big(\lambda_q, q^{-s}\big)^{t_q+1}} \prod_{p \mid M} \frac{p^{-s t_p}}{\mathrm{L}_p\big(\lambda_p, p^{-s}\big)^{t_p+1}} \pmod{\ell}. \qquad \blacksquare$$

In Sections 3.4 and 3.7, we combine Corollary 3.3 with the formula for the L-series of Shimura lifts in (1.6). To this end, we introduce notation for the relevant products of inverse $\mathrm{L}_p$-polynomials. Given eigenvalues $\lambda = (\lambda_p)_{p \mid M}$ modulo $\ell$ and nonnegative integers $t = (t_p)_{p \mid M}$, we let

$$c_M(\lambda, t;\, m) := c\Bigg( \prod_{p \mid M} \frac{p^{-s t_p}\big(1 - \chi \chi_{-4}^{2k-r} \chi_D(p) p^{k - \frac{3}{2} - s}\big)}{\big(1 - \lambda_p p^{-s} + \chi^2(p) p^{2k-2-2s}\big)^{t_p+1}};\, m \Bigg) \tag{3.10}$$

be the coefficient of $m^{-s}$ in the given Dirichlet series. Note that $k$, $\chi$, and $D$ are suppressed from our notation and given by context. Observe also that these coefficients are trivial if $m$ does not divide some power of $M$.

**3.4 $\mathrm{U}_p$-congruences** The next proposition is concerned with $\mathrm{U}_p$ congruences. It implies that maximal Ramanujan-type congruences have $M_p \nmid \beta$ for every $p \mid M$.

**Proposition 3.5.** *Fix a cusp form $f$ with a Ramanujan-type congruence as in* (3.2). *For any prime $p \mid M$ with $M_p \mid \beta$, we have the Ramanujan-type congruence*

$$\forall n \in \mathbb{Z} :\; c(f;\, M_p^{\#} n + \beta) \equiv 0 \pmod{\ell}.$$

*Proof.* In the following argument, we can restrict to the case of the theta-character after replacing $f$ with $\mathrm{V}_{24} f$, if $f$ is associated with the eta-character. We further employ Lemma 2.2 in conjunction with Proposition 2.1 and replace $f$ by

$$\sum_{\substack{n \in \mathbb{Z} \\ n \,\square\, \beta \pmod{M_p^{\#}}}} c(f;\, n) e(n\tau).$$

This allows us to assume that $M = M_p$, and we have to show that $f \equiv 0 \pmod{\ell}$.

Since $M \mid \beta$, we have $M\mathbb{Z} + \beta = M\mathbb{Z}$. In other words, we have $c(f; Mn) \equiv 0$ for all positive integers $n$. Recall from Section 3.1 that $M = M_{\mathrm{fd}} M_{\mathrm{s}}^2$ where $M_{\mathrm{fd}} \in \{1, p\}$ and $M_{\mathrm{s}}$ is a suitable power of $p$. Given any integer $m$ and any fundamental discriminant $D$ with $\mathrm{sgn}(D) = (-1)^{k-\frac{1}{2}}$, we have the following congruences:

$$\begin{aligned} c(f; |D| M_{\mathrm{s}}^2 m^2) &\equiv 0 \pmod{\ell}, \quad \text{if } p \mid D \text{, or if } p \nmid D \text{ and } p \nmid M_{\mathrm{fd}}; \\ c(f; |D| M_{\mathrm{s}}^2 p^2 m^2) &\equiv 0 \pmod{\ell}, \quad \text{if } p \nmid D \text{ and } p \mid M_{\mathrm{fd}}. \end{aligned}$$

To simplify the argument, we will employ the possibly weaker congruence for the coefficients of index $|D| M_{\mathrm{s}}^2 p^2 m^2$.

Recall the description of the L-series of the $D$-th Shimura lift in (1.6), where $D$ is a fundamental discriminant with $\mathrm{sgn}(D) = (-1)^{k-\frac{1}{2}}$. The congruences for the Fourier coefficients of $f$ that we have found show the coefficients of $(M_{\mathrm{s}} pm)^{-s}$ of the L-series on the right hand side of (1.6) are congruent to zero modulo $\ell$.

We insert the decomposition of the L-series of the $D$-th Shimura lift of $f$ in Corollary 3.3 into the left hand side of (1.6). We thus obtain congruences for the coefficients

$$c\Bigg( \sum_{\lambda} \sum_{t=0}^{d_\lambda} \frac{\mathrm{L}^p\big(\mathrm{Sh}_D(f_{\lambda,t}),\, s\big)}{\mathrm{L}^p\big(\chi \chi_{-4}^{2k-r} \chi_D,\, s + \frac{3}{2} - k\big)} \, \frac{p^{-s t_p} \big(1 - \chi \chi_4^{2k-r} \chi_D(p) p^{k-\frac{3}{2}-s}\big)}{\mathrm{L}_p\big(\lambda_p, p^{-s}\big)^{t_p+1}} ;\, M_{\mathrm{s}} pm \Bigg)$$

that vanish modulo $\ell$.

Now assume by contraposition that we have $f \not\equiv 0 \pmod{\ell}$. Then we can fix some $D$ such that $\mathrm{Sh}_D(f) \not\equiv 0 \pmod{\ell}$. The generalized eigenspaces modulo $\ell$ of $\mathrm{T}_p^{\mathrm{cl}}$ intersect trivially, as they are the direct summands in the primary decomposition of the corresponding $\mathrm{F}_{K,\ell}[\mathrm{T}_p^{\mathrm{cl}}]$-module. Since further $m$ is arbitrary and the denominator $\mathrm{L}^p(\chi \chi_D, s + \frac{3}{2} - k)$ is common to all terms in the double sum over $\lambda$ and $t$, we can invoke Lemma 2.6 to conclude congruence modulo $\ell$ for the coefficients

$$c\Bigg( \sum_{t=0}^{d_\lambda} \frac{\mathrm{L}^p\big(\mathrm{Sh}_D(f_{\lambda,t}),\, s\big)}{\mathrm{L}^p\big(\chi \chi_{-4}^{2k-r} \chi_D,\, s + \frac{3}{2} - k\big)} \, \frac{p^{-s t_p} \big(1 - \chi \chi_{-4}^{2k-r} \chi_D(p) p^{k-\frac{3}{2}-s}\big)}{\mathrm{L}_p\big(\lambda_p, p^{-s}\big)^{t_p+1}} ;\, M_{\mathrm{s}} pm \Bigg) \equiv 0 \pmod{\ell}.$$

Observe that $f_{\lambda,0}$ is linearly independent modulo $\ell$ of the span of all $f_{\lambda,t}$ with $t \ge 1$. Applying Lemma 2.6 once more, we find a congruence modulo $\ell$ for the coefficients

$$c\Bigg( \frac{\mathrm{L}^p\big(\mathrm{Sh}_D(f_{\lambda,0}),\, s\big)}{\mathrm{L}^p\big(\chi \chi_{-4}^{2k-r} \chi_D,\, s + \frac{3}{2} - k\big)} \, \frac{1 - \chi \chi_{-4}^{2k-r} \chi_D(p) p^{k-\frac{3}{2}-s}}{\mathrm{L}_p\big(\lambda_p, p^{-s}\big)} ;\, M_{\mathrm{s}} pm \Bigg) \equiv 0 \pmod{\ell}.$$

Specializing this congruences to the case when $m$ is a power of $p$, we find that

$$\frac{1-\chi\chi_{-4}^{2k-r}\chi_D(p)p^{k-\frac{3}{2}-s}}{\mathrm{L}_p\big(\lambda_p, p^{-s}\big)} \;\equiv\; R(p^{-s}) \;\; (\mathrm{mod}\,\ell)$$

for some polynomial $R(X)$ over $\mathcal{O}_{K,\ell}$ with constant coefficient 1.

Since $p \nmid N$, we can factor $L_p(\lambda_p, X) = (1-\alpha_p X)(1-\beta_p X)$ over $\mathbb{C}$. In the remainder of the proof, we can and will assume that $\alpha_p, \beta_p \in \mathcal{O}_{K,\ell}$ by extending $K$ if needed. Since $p$ is co-prime to $\ell$, we conclude that $\alpha_p, \beta_p \not\equiv 0 \;(\mathrm{mod}\,\ell)$. Rephrasing the previous congruence, we have

$$1-\chi\chi_{-4}^{2k-r}\chi_D(p)p^{k-\frac{3}{2}}X \;\equiv\; R(X)\,(1-\alpha_p X)(1-\beta_p X) \;\; (\mathrm{mod}\,\ell).$$

When comparing degrees, we obtain a contradiction. We therefore must have the desired congruence $f \equiv 0 \;(\mathrm{mod}\,\ell)$. ∎

**3.5 Fundamental discriminants and square-classes** Consider a Ramanujan-type congruence as described in (3.2), and recall the definition of $M_0$, $\beta_0$, and $M_1$ from Section 3.1. We define the set of fundamental discriminants related to the arithmetic progression $M\mathbb{Z}+\beta$ as

$$\begin{aligned}\mathrm{D}_k\big(M\mathbb{Z}+\beta\big) \;:=\; \big\{&D \text{ fund. disc.} : \mathrm{sgn}(D) = (-1)^{k-\frac{1}{2}},\\ &\exists m \in \mathbb{Z}, n_0 \in \mathbb{Z}, n_0 \,\square\, \beta_0 \;(\mathrm{mod}\, M_0) : |D|m^2 = M_1 n_0\big\}.\end{aligned} \tag{3.11}$$

The next lemma provides an explicit description of $\mathrm{D}_k(M\mathbb{Z}+\beta)$ and relates the integers that occur in the multiplicative congruence (3.5) with the fundamental discriminants in (3.11). We use the notation $M_{\mathrm{fd}}$ and $M_{\mathrm{s}}$ defined in Section 3.1.

**Lemma 3.6.** *Fix a weight-k cusp form $f$ with a maximal Ramanujan-type congruence on $M\mathbb{Z}+\beta$ as in* (3.2)*. Then we have*

$$\mathrm{D}_k(M\mathbb{Z}+\beta) \;=\; \big\{D \textit{ fund. disc.} : \mathrm{sgn}(D) = (-1)^{k-\frac{1}{2}},\, M_{\mathrm{fd}} \,|\, D,\, \big(|D|/M_{\mathrm{fd}}\big) \,\square\, \beta_0 \;(\mathrm{mod}\, M)\big\}$$

*and*

$$\begin{aligned}&\big\{M_1 n_0 \,:\, n_0 \in \mathbb{Z}, n_0 > 0, n_0 \,\square\, \beta_0 \;(\mathrm{mod}\, M)\big\}\\ &\qquad = \big\{|D| M_{\mathrm{s}}^2 m^2 \,:\, D \in \mathrm{D}_k(M\mathbb{Z}+\beta),\, m \in \mathbb{Z},\, \gcd(m, M) = 1\big\}.\end{aligned}$$

*Proof.* In the case of the theta-character $M$ is odd since it is co-prime to $N$, in the case of the eta-character it is odd, since we assume that the Ramanujan-type congruences (3.2) is maximal. We find that in particular $M_{\mathrm{fd}} \,|\, M$ is odd. The maximality of the Ramanujan-type congruences and Proposition 3.5 also imply that every prime $p$ that divides $M$ equally divides $M_0$ but not $\beta_0$. This already implies the given description of $\mathrm{D}_k(M\mathbb{Z} + \beta)$.

Consider the second equality in the lemma. We directly verify that the right hand side is contained in the left hand side. Conversely, consider any element of the left hand side. We can write it as $|D|m'^2$ for a fundamental discriminant $D \in \mathrm{D}_k(M\mathbb{Z} + \beta)$ and a positive integer $m'$. Specifically, we have $|D|m'^2 = M_{\mathrm{fd}} M_{\mathrm{s}}^2 n_0$ for some positive integer $n_0$ with $n_0 \square \beta_0 \pmod{M_0}$. Since both $M_{\mathrm{fd}} M_{\mathrm{s}}^2$ and $\gcd(\beta_0, M_0)$ are odd, we have $M_{\mathrm{fd}} \,|\, D$, $M_{\mathrm{s}} \,|\, m'$, and $\gcd(n_0, M_0) = 1$. Further, since $n_0 \square n_0 n^2 \pmod{M_0}$ for every integer $n$ with $\gcd(n, M_0) = 1$, we have $m' = M_{\mathrm{s}} m$ for some integer $m$ that satisfies $\gcd(m, M_0) = 1$. We conclude that $\gcd(m, M) = 1$ by using once more that every prime $p \,|\, M$ also divides $M_0$. ∎

**3.6 Non-zero Shimura lifts** We next describe a dichotomy between Ramanujan-type congruences that can be related to congruences for modular forms of integral weight via a Shimura lift and Ramanujan-type congruences involving fundamental discriminants.

Note that here and in the remainder we shall apply the Shimura lift modulo $\ell$ in (1.7) and (1.8) to modular forms modulo $\ell$. In particular, we can apply it to cusp forms of weight $k$ less than $\frac{5}{2}$. While the resulting modular form might have weight larger than $2k - 1$, its weight is congruent to $2k - 1$ modulo $\ell_{\mathbb{Z}} - 1$, where $\ell_{\mathbb{Z}}\mathbb{Z} = \ell \cap \mathbb{Z}$.

In the statement of the next proposition, we use the notation defined in Section 3.5.

**Proposition 3.7.** *Fix a cusp form $f$ with a Ramanujan-type congruence as in (3.2). Either there is at least one discriminant $D \in \mathrm{D}_k(M\mathbb{Z} + \beta)$ such that $\mathrm{Sh}_D(f) \not\equiv 0 \pmod{\ell}$ or we have the congruences*

$$\forall\, n_0 \in \mathbb{Z}, n_0 \square \beta_0 \pmod{M_0} : c(f; M_{\mathrm{fd}} n_0) \equiv 0 \pmod{\ell}. \tag{3.12}$$

*Remark 3.8.* The congruence (3.12) subsumes the Ramanujan-type congruence on the arithmetic progression $M_{\mathrm{fd}}(M_0\mathbb{Z} + \beta_0)$. The crucial point in light of scarcity of congruences investigated in previous work [2] is that $M_{\mathrm{fd}} M_0$ is cube-free.

*Proof.* Suppose that we have $\mathrm{Sh}_D(f) \equiv 0 \pmod{\ell}$ for all $D \in D_k(M\mathbb{Z} + \beta)$. Then the coefficient formula in (1.6) for the Shimura lift implies that $c(f; |D|m^2) \equiv 0 \pmod{\ell}$ for all $D \in D_k(M\mathbb{Z} + \beta)$ and all $m \in \mathbb{Z}$.

Fix some $n_0 \in \mathbb{Z}$ with $n_0 \square \beta_0 \pmod{M_0}$ as in (3.12). Along the lines of the proof of Lemma 3.6, we see that $M_{\mathrm{fd}} n_0 = |D| m^2$ for some $D \in D_k(M\mathbb{Z} + \beta)$ and some $m \in \mathbb{Z}$ with $\gcd(m, M_0) = 1$. We conclude that $c(f; M_{\mathrm{fd}} n_0) = c(f; |D| m^2) \equiv 0 \pmod{\ell}$, confirming the claim. ■

**3.7 Congruences for coefficients of L-series** In this section, we characterize Ramanujan-type congruences for modular forms by congruences for the Dirichlet series coefficients in (3.10).

**Lemma 3.9.** *Fix a cusp form $f$ with a Ramanujan-type congruence as in (3.2), which we assume is maximal. Given a fundamental discriminant $D \in \mathrm{D}_k(M\mathbb{Z} + \beta)$, a tuple of eigenvalues $\lambda = (\lambda_p)_{p \mid M}$ modulo $\ell$, and an $\mathrm{F}_{K,\ell}$-linear map $\phi \in \mathrm{Sh}_D(F_\lambda)^\vee$, we have the congruence*

$$\sum_{t=0}^{d_\lambda} \phi\big(\mathrm{Sh}_D(f_{\lambda,t})\big) c_M(\lambda, t; M_{\mathrm{s}}) \equiv 0 \pmod{\ell}. \tag{3.13}$$

*Conversely, consider a modular form as in (3.1). Given a positive integer $M$ and an integer $\beta$ with $\gcd(M, \ell N) = 1$, $M \mid M_{\mathrm{sf}} \beta$, and $M_p \nmid \beta$ for every prime $p \mid M$, assume that the congruence in (3.13) holds for every $D \in D(M\mathbb{Z} + \beta)$, $\lambda$, and $\phi \in \mathrm{Sh}_D(F_\lambda)^\vee$. Then we have a Ramanujan-type congruence as in (3.2).*

*Proof.* Similar to the proof of Proposition 3.5, we combine the equation in (1.6) for the L-series of Shimura lifts with Corollary 3.3 to find that

$$\begin{aligned} &\sum_{n=1}^{\infty} c\big(f; |D| n^2\big) n^{-s} \\ &\equiv \sum_{\lambda} \sum_{t=0}^{d_\lambda} \frac{\mathrm{L}^M\big(\mathrm{Sh}_D(f_{\lambda,t}), s\big)}{\mathrm{L}^M\big(\chi\chi_{-4}^{2k-r}\chi_D, s + \frac{3}{2} - k\big)} \prod_{p \mid M} \frac{p^{-s t_p}\big(1 - \chi\chi_{-4}^{2k-r}\chi_D(p) p^{k-\frac{3}{2}-s}\big)}{L_p\big(\lambda_p, p^{-s}\big)^{t_p+1}} \pmod{\ell}. \end{aligned}$$

The multiplicative congruence in (3.5) and Lemma 3.6 imply that, for all integers $m$ with $\gcd(m, M) = 1$, we have $c(f; |D| M_{\mathrm{s}}^2 m^2) \equiv 0 \pmod{\ell}$. We infer from this that

$$\sum_{\lambda} \sum_{t=0}^{d_\lambda} \frac{\mathrm{L}^M\big(\mathrm{Sh}_D(f_{\lambda,t}), s\big)}{\mathrm{L}^M\big(\chi\chi_{-4}^{2k-r}\chi_D, s + \frac{3}{2} - k\big)} c_M(\lambda, t; M_{\mathrm{s}}) \equiv 0 \pmod{\ell}.$$

We can multiply this congruence with $\mathrm{L}^M(\chi\chi_{-4}^{2k-r}\chi_D, s + \frac{3}{2} - k)$, and thus cancel the denominator on the left hand side.

Lemma 2.6 allows us to replace the Dirichlet series $\mathrm{L}^M(\mathrm{Sh}_D(f_{\lambda,t}), s)$ by $\mathrm{Sh}_D(f_{\lambda,t})$. We find that

$$\sum_{\lambda}\sum_{t=0}^{d_\lambda} \mathrm{Sh}_D(f_{\lambda,t})\, c_M(\lambda, t; M_{\mathrm{s}}) \;\equiv\; 0 \;\; (\mathrm{mod}\, \ell).$$

Since the generalized eigenspaces for the Hecke operators $\mathrm{T}^{\mathrm{cl}}_p$, $p \,|\, M$ intersect pairwise trivially, each term in the sum over $\lambda$ is congruent to zero modulo $\ell$. The first part of the statement follows by applying $\phi$ to the resulting congruence.

The converse follows by applying the argument in reverse order, where the conditions on $M$ are needed to obtain a characterization as in Lemma 3.6. ∎

As an immediate consequence of Lemma 3.9, we can restrict our attention to single tuples of eigenvalues $\lambda = (\lambda_p)_{p\,|\,M}$, specified in Section 3.2.

**Corollary 3.10.** *Fix a cusp form $f$ with a Ramanujan-type congruence as in* (3.2)*. Then we have*

$$\forall \lambda = (\lambda_p)_{p\,|\,M} \;(\mathrm{mod}\, \ell)\; \forall n \in \mathbb{Z} : c(f_\lambda;\, Mn + \beta) \equiv 0 \;\; (\mathrm{mod}\, \ell).$$

**3.8 Hecke modules** We investigate the sets of modular forms with a given Ramanujan-type congruence. Adjusting to the setup in (3.1) and (3.2), we consider a positive integer $M$ and an integer $\beta$ with $\gcd(M, \ell N) = 1$. To ease notation, we set $\widetilde{\chi} = \chi \chi_{\theta+}^r$, $r$ odd, or $\widetilde{\chi} = \chi\chi_\eta^r$, $\gcd(r, 24) = 1$. We define the $\mathcal{O}_{K,\ell}$-module

$$\mathrm{R}_\ell\big(M\mathbb{Z} + \beta, \mathrm{S}_k(\widetilde{\chi})\big) \;:=\; \big\{ f \in \mathrm{S}_k(\widetilde{\chi}; \mathcal{O}_{K,\ell}) : \forall n \in \mathbb{Z} : c(f; Mn + \beta) \equiv 0 \;(\mathrm{mod}\, \ell) \big\}. \qquad (3.14)$$

As in the case of Ramanujan-type congruences for individual modular forms, there is no loss in generality when employing Proposition 2.1 to assume that $M \,|\, M_{\mathrm{sf}}\beta$. We can also assume that $M_p \nmid \beta$ for every prime $p \,|\, M$ by applying Proposition 3.5.

The next proposition states that the $\mathcal{O}_{K,\ell}$-module in (3.14) is a Hecke-module away from $\gcd(M, \ell N)$.

**Proposition 3.11.** *Consider the set of cusp forms with a Ramanujan-type congruence as in* (3.14)*. Given a prime $p$ with* $\gcd(p, M, \ell N) = 1$*, we have*

$$\mathrm{R}_\ell\big(M\mathbb{Z} + \beta, \mathrm{S}_k(\widetilde{\chi})\big) \big|\, \mathrm{T}^{\mathrm{cl}}_p \;\subseteq\; \mathrm{R}_\ell\big(M\mathbb{Z} + \beta, \mathrm{S}_k(\widetilde{\chi})\big). \qquad (3.15)$$

*Proof.* For consistency with other proofs, we write $q$ instead of $p$ in this proof. For those primes $q$ with $q \,|\, N$ or $\ell \,|\, q$, the action of the Hecke operator modulo $\ell$ coincides with the one of $\mathrm{U}_{q^2}$, and the result follows directly, since multiplication by $q^2$, $q \nmid M$, preserves square-classes modulo $M$. In the remainder of the proof, we can therefore assume that $\gcd(q, \ell N) = 1$.

Given any $f$ as in (3.2), we inspect the Ramanujan-type congruences of $\mathrm{Sh}_D(f)$ for all $D \in \mathrm{D}_k(M\mathbb{Z} + \beta)$ to show that $f|\mathrm{T}^{\mathrm{cl}}_q$ has a Ramanujan-type congruence on $M\mathbb{Z} + \beta$. If we have $\mathrm{Sh}_D(f) \equiv 0 \pmod{\ell}$, we find that

$$\mathrm{Sh}_D\big(f|\mathrm{T}^{\mathrm{cl}}_q\big) \;=\; \mathrm{Sh}_D(f)\big|\mathrm{T}^{\mathrm{cl}}_q \;\equiv\; 0 \pmod{\ell}.$$

In particular, we have $c(f|\mathrm{T}^{\mathrm{cl}}_q; Mn{+}\beta) \equiv 0 \pmod{\ell}$ for all $n \in \mathbb{Z}$ such that there is $m \in \mathbb{Z}$ with $Mn + \beta = |D| m^2$. We can therefore assume that $\mathrm{Sh}_D(f)$ is not congruent to zero modulo $\ell$.

By Corollary 3.10 it suffices to treat the case $f = f_\lambda$ for a fixed tuple of eigenvalues $\lambda = (\lambda_p)_{p \,|\, M}$ modulo $\ell$. The case of $\gcd(q, \ell M N) = 1$ follows immediately by applying the converse direction in Lemma 3.9.

Assume that $q \,|\, M$. We use induction on the Jordan-Hölder lengths of the module generated by the action of $\mathrm{T}^{\mathrm{cl}}_p$, $p \,|\, M$, on $\mathrm{F}_{K,\ell}\, f$. In the base case of the induction $f$ is an eigenform modulo $\ell$ for all these $\mathrm{T}^{\mathrm{cl}}_p$, and the statement follows trivially. We set $g := f|(\mathrm{T}^{\mathrm{cl}}_q - \lambda_q)$. By induction, we can assume that the Hecke operators $\mathrm{T}^{\mathrm{cl}}_p$, $p \,|\, M$, preserve any Ramanujan-type congruence satisfied by $g$. In particular, it suffices to show that $g$ has a Ramanujan-type congruence modulo $\ell$ on $M\mathbb{Z} + \beta$.

For tuples of integers $t = (t_p)_{p \,|\, M}$ define the modular forms as in (3.6) by

$$g_t \;:=\; g\Big|\prod_{p \,|\, M} \big(\mathrm{T}^{\mathrm{cl}}_p - \lambda_p\big)^{t_p},$$

and let $G$ be the span over $\mathrm{F}_{K,\ell}$ of all $g_t \pmod{\ell}$. The desired Ramanujan-type congruence for $g$ follows via Lemma 3.9, if we show that for every $\psi \in G^\vee$, we have

$$\sum_{t=0}^{d_\lambda} \psi\big(\mathrm{Sh}_D(g_t)\big) c_M(\lambda, t;\, M_{\mathrm{s}}) \;\equiv\; 0 \pmod{\ell}.$$

For clarity, we remark that in this sum some terms are trivially zero, since the $d_{\lambda,p}$-th power of $\mathrm{T}^{\mathrm{cl}}_q - \lambda$ annihilates $g$ modulo $\ell$. Let $\phi := \psi \circ (\,\cdot\,|(\mathrm{T}^{\mathrm{cl}}_q - \lambda)) \in F^\vee_\lambda$. Then Lemma 3.9 asserts that

$$\sum_{t=0}^{d_\lambda} \phi\big(\mathrm{Sh}_D(f_{\lambda,t})\big) c_M(\lambda, t;\, M_{\mathrm{s}}) \;\equiv\; 0 \pmod{\ell}.$$

Since $\phi(\mathrm{Sh}_D(f_{\lambda,t})) = \psi(\mathrm{Sh}_D(g_t))$, the statement follows. ∎

**3.9 Ramanujan-type congruences for eigenforms** In this section, we clarify the relation between Ramanujan-type congruences and congruences satisfied by Hecke eigenvalues modulo $\ell$.

**Proposition 3.12.** *Fix a cusp form $f$ as in (3.1), an odd prime $p$, $\gcd(p, \ell N) = 1$, and $\epsilon \in \{0, \pm 1\}$. In the case of the eta-character assume that $p \ne 3$. Suppose that there is a fundamental discriminant $D$ satisfying $\mathrm{sgn}(D) = (-1)^{k-\frac{1}{2}}$ and $\chi_D(p) = \epsilon$ such that we have $\mathrm{Sh}_D(f) \not\equiv 0 \pmod{\ell}$. Assume further that for all such $D$ we have the congruence*

$$\mathrm{Sh}_D(f) \big| \mathrm{T}_p^{\mathrm{cl}} \equiv \lambda_p \mathrm{Sh}_D(f) \pmod{\ell}, \quad \lambda_p \in \mathcal{O}_{K,\ell}.$$

*(1) Consider the case that $\epsilon = 0$. If $\lambda_p^2 \equiv 4\chi^2(p) p^{2k-2} \pmod{\ell}$, we have the Ramanujan-type congruence with gaps*

$$\forall n \in \mathbb{Z} \setminus p\mathbb{Z} : c(f;\, p^{2m+1} n) \equiv 0 \pmod{\ell} \tag{3.16}$$

*for a positive integer $m$ if and only if $m \equiv -1 \pmod{\ell}$. If $\lambda_p^2 \not\equiv 4\chi^2(p) p^{2k-2} \pmod{\ell}$ and the $\mathrm{L}$-polynomial*

$$1 - \lambda_p X + \chi^2(p) p^{2k-2} X^2 \equiv (1 - \alpha_p X)(1 - \beta_p X) \pmod{\ell}$$

*factors over $\mathrm{F}_{K,\ell}$, then (3.16) holds if and only if $\alpha_p^{m+1} \equiv \beta_p^{m+1} \pmod{\ell}$. If $m = 1$, this is equivalent to $\lambda_p \equiv 0 \pmod{\ell}$.*

*(2) Consider the case that $\epsilon \in \{\pm 1\}$. If $\lambda_p^2 \equiv 4\chi^2(p) p^{2k-2} \pmod{\ell}$, we have the multiplicative congruence*

$$\forall n \in \mathbb{Z}, \left( \frac{(-1)^{k-\frac{1}{2}} n}{p} \right) = \epsilon : c(f;\, p^{2m} n) \equiv 0 \pmod{\ell} \tag{3.17}$$

*for a positive integer $m$ if and only if*

$$m \not\equiv -1 \pmod{\ell}, \quad \lambda_p \equiv \frac{2m \chi_{-4}^{2k-r} \chi(p) \epsilon p^{k-\frac{3}{2}}}{m+1} \pmod{\ell}, \quad \text{and} \quad p \equiv \frac{m^2}{(m+1)^2} \pmod{\ell}. \tag{3.18}$$

*If $\lambda_p^2 \not\equiv 4\chi^2(p) p^{2k-2} \pmod{\ell}$ and the $\mathrm{L}$-polynomial factors as before, then (3.16) holds if and only if*

$$\frac{\alpha_p^m}{\beta_p^m} \equiv \frac{\beta_p - \chi_{-4}^{2k-r} \chi(p) \epsilon p^{k-\frac{3}{2}}}{\alpha_p - \chi_{-4}^{2k-r} \chi(p) \epsilon p^{k-\frac{3}{2}}} \pmod{\ell}. \tag{3.19}$$

*If $m = 1$, this congruence is equivalent to $\lambda_p \equiv \chi_{-4}^{2k-r} \chi(p) \epsilon p^{k-\frac{3}{2}} \pmod{\ell}$.*

*Remark 3.13.* The characterization that we provide for $p \mid D$ is well-suited to both establish and disprove Ramanujan-type congruences by for instance inspecting the multiplicative order of $\alpha_p / \beta_p$. In particular, in these cases there is always some $m$ for which there is a Ramanujan-type congruence. The situation is different if $p \nmid D$. There might be modular forms without any multiplicative congruence (3.17). The characterization that we provide is still very suitable to rule them out. At the end of Example 3.19, we illustrate this process.

*Proof of Proposition 3.12.* By Lemma 3.9 and using notation defined in (3.10), we have to characterize the congruences $c_p(\lambda, 0; p^m) \equiv 0 \pmod{\ell}$ for fundamental discriminants $D$, which are suppressed from our notation, with $\chi_D(p) = 0$ in the first case and $\chi_D(p) = \epsilon$ in the second one. Observe that the first case corresponds exactly to the case of integral weight congruences of Shimura lifts. We can therefore adopt the result from Proposition 4.6 in [18], where the weight is $2k-1$ and the Dirichlet character is $\chi^2$. The statement for $m = 1$, follows since $\alpha_p \equiv \pm\beta_p \pmod{\ell}$, but $\alpha_p \not\equiv \beta_p \pmod{\ell}$ by the assumption on $\lambda_p^2$. Then we have $\lambda_p \equiv \alpha_p + \beta_p \equiv 0 \pmod{\ell}$.

We inspect the case of $p \nmid D$ and need to determine the coefficients $c_p(\lambda, 0; p^m)$ of

$$\frac{1 - \chi_{-4}^{2k-r} \chi(p) \epsilon p^{k-\frac{3}{2}} X}{1 - \lambda_p X + \chi^2(p) p^{2k-2} X^2}.$$

By extending $K$ we can assume that the L-polynomial factors as in the statement. Consider the case that $\lambda_p^2 \equiv 4\chi^2(p) p^{2k-2} \pmod{\ell}$, that is, $\alpha_p \equiv \beta_p \pmod{\ell}$. The $m$-th coefficient of the inverse of $1 - \lambda_p X + \chi^2(p) p^{2k-2} X^2$ is congruent to $(m+1)\alpha_p^m$. Since we have $m \ge 1$, we obtain that

$$c_p(\lambda, 0; p^m) \;\equiv\; \alpha_p^{m-1} \big( (m+1)\alpha_p - m \chi_{-4}^{2k-r} \chi(p) \epsilon p^{k-\frac{3}{2}} \big) \pmod{\ell}.$$

Since $\alpha_p^{m-1} \not\equiv 0 \pmod{\ell}$, we conclude that the multiplicative congruence (3.17) is characterized by

$$m \not\equiv -1 \pmod{\ell} \quad \text{and} \quad \lambda_p \equiv 2\alpha_p \equiv \frac{2m \chi_{-4}^{2k-r} \chi(p) \epsilon p^{k-\frac{3}{2}}}{m+1} \pmod{\ell}.$$

From $\alpha_p^2 \equiv \chi^2(p) p^{2k-2} \pmod{\ell}$, we infer that $p \equiv m^2 / (m+1)^2 \pmod{\ell}$.

Consider the case that $\lambda_p^2 \not\equiv 4\chi^2(p) p^{2k-2} \pmod{\ell}$, that is, $\alpha_p \not\equiv \beta_p \pmod{\ell}$. Then a calculation shows that the $m$-th coefficient of the inverse of $1 - \lambda_p X + \chi^2(p) p^{2k-2} X^2$ is congruent to $(\alpha_p^{m+1} - \beta_p^{m+1}) / (\alpha_p - \beta_p)$. We obtain the congruence

$$c_p(\lambda, 0; p^m) \equiv \frac{\alpha_p^m (\alpha_p - \chi_{-4}^{2k-r} \chi(p) \epsilon p^{k-\frac{3}{2}}) - \beta_p^m (\beta_p - \chi_{-4}^{2k-r} \chi(p) \epsilon p^{k-\frac{3}{2}})}{\alpha_p - \beta_p} \pmod{\ell}.$$

Since $\alpha_p \not\equiv \beta_p \pmod{\ell}$, we conclude that the multiplicative congruence (3.17) is characterized by

$$\alpha_p^m \big(\alpha_p - \chi_{-4}^{2k-r} \chi(p) \epsilon p^{k-\frac{3}{2}}\big) \equiv \beta_p^m \big(\beta_p - \chi_{-4}^{2k-r} \chi(p) \epsilon p^{k-\frac{3}{2}}\big) \pmod{\ell},$$

as stated in the proposition.

To obtain the characterization in the case of $m = 1$, it is most helpful to calculate directly the first coefficient of

$$\frac{1 - \chi_{-4}^{2k-r} \chi(p) \epsilon p^{k-\frac{3}{2}} X}{1 - \lambda_p X + \chi^2(p) p^{2k-2} X^2} = 1 + \big(\lambda_p - \chi_{-4}^{2k-r} \chi(p) \epsilon p^{k-\frac{3}{2}}\big) X + \cdots. \qquad \blacksquare$$

**Corollary 3.14.** *Consider the situation of Proposition 3.12. If $m = 1$, we cannot have the multiplicative congruence in (3.17) for both $\epsilon = +1$ and $\epsilon = -1$.*

*Proof.* We use the same notation as in Proposition 3.12. Regardless of whether we have $\alpha_p \equiv \beta_p \pmod{\ell}$ or $\alpha_p \not\equiv \beta_p \pmod{\ell}$, the proposition yields the congruence

$$\chi_{-4}^{2k-r} \chi(p) p^{k-\frac{3}{2}} \equiv -\chi_{-4}^{2k-r} \chi(p) p^{k-\frac{3}{2}} \pmod{\ell},$$

which is impossible, since $\ell \nmid 2$ and $p \not\equiv 0 \pmod{\ell}$. $\blacksquare$

**3.10 Ramanujan-type congruences with gap** In Section 2.2, we have seen that if the square of a prime exactly divides $M$ for a maximal Ramanujan-type congruence modulo $\ell$ on $M\mathbb{Z} + \beta$, we find congruences on two square-classes modulo $M$. In this section, we extend the result of Section 2.2 to all squares that exactly divide $M$.

**Proposition 3.15.** *Fix a cusp form $f$ with a maximal Ramanujan-type congruence as in (3.2). Assume that $M_p \ne p^2$ is a square for a prime $p \mid M$. Then we have the Ramanujan-type congruences*

$$\forall \beta' \in \mathbb{Z}, \gcd(\beta', M_p) = \gcd(\beta, M_p),\ \beta' \,\square\, \beta \pmod{M_p^\#} :$$
$$\forall n \in \mathbb{Z} : c(f; Mn + \beta') \equiv 0 \pmod{\ell}.$$

*Proof.* Corollary 3.10 allows us to assume that $f$ is a generalized eigenform modulo $\ell$ with respect to $\mathrm{T}_p^{\mathrm{cl}}$. Observe that we have $p \,|\, D$ for all $D \in \mathrm{D}_k(M\mathbb{Z} + \beta)$, since $M_p$ is a square that exactly divides $p\beta$. We therefore have $\chi_D(p) = 0$ in the defining equation (3.10) of the coefficients $c_p(\lambda, t; M_{\mathrm{s}})$ that appear in Lemma 3.9. Using the converse direction of the same lemma, our proposition follows. ∎

**3.11 Prime-power Ramanujan-type congruences** Many cases of Ramanujan-type congruences in the literature, say for the partition function function, arise from Ramanujan-type congruences of cusp forms modulo $\ell$, for instance $\mathrm{U}_\ell \eta^{-1} \pmod{\ell}$, on prime powers $M$. In this section, we show that for Hecke eigenforms modulo $\ell$ as opposed to generalized Hecke eigenforms maximal Ramanujan-type congruences satisfy $M = M_p$. We also give examples of Ramanujan-type congruences that do not arise in this way, illustrating that neither proper generalized Hecke eigenforms modulo $\ell$ nor linear combinations of Hecke eigenforms modulo $\ell$ satisfy the conclusion of the next proposition.

**Proposition 3.16.** *Fix a cusp form $f$ with a Ramanujan-type congruence as in (3.2). Assume that $f$ is an eigenform modulo $\ell$ for all $\mathrm{T}_p^{\mathrm{cl}}$, $p \,|\, M$. Then there is a prime $p \,|\, M$ such that*

$$\forall n \in \mathbb{Z} : c(f; M_p n + \beta) \equiv 0 \pmod{\ell}.$$

*Proof.* We can assume that $f \not\equiv 0 \pmod{\ell}$. For every $D \in \mathrm{D}_k(M\mathbb{Z} + \beta)$, we employ the factorization in Lemma 3.9. Since $f$ is a Hecke eigenform modulo $\ell$ for $\mathrm{T}_p^{\mathrm{cl}}$, $p \,|\, M$, the coefficients $c_M(\lambda, 0; m)$ that appear in Lemma 3.9 and are defined in (3.10) are multiplicative in $m$. In other words, we have

$$c_M(\lambda, 0; m) \equiv \prod_{p \,|\, M} c_M(\lambda, 0; m_p) \pmod{\ell},$$

where $m = \prod m_p$ is a factorization into powers $m_p$ of the primes $p \,|\, M$.

We consider $m = M_{\mathrm{s}}$, for which we have $c_M(\lambda, 0; M_{\mathrm{s}}) \equiv 0 \pmod{\ell}$ by Lemma 3.9. Fix some $p \,|\, M$ such that $c_p(\lambda_p, 0; m_p) \equiv 0 \pmod{\ell}$. Applying the converse direction in Lemma 3.9 for every $D \in \mathrm{D}_k(M_p\mathbb{Z} + \beta)$, we confirm the statement. ∎

**3.12 Examples** We finish the section by several examples of cusp forms that have a Ramanujan-type congruence on $p_1^{m_1} p_2^{m_2} \mathbb{Z} + \beta$ that neither extends to $p_1^{m_1} \mathbb{Z} + \beta$ nor to $p_2^{m_2} \mathbb{Z} + \beta$.

**Example 3.17.** We consider the sum of two eigenforms to find our first example. Specifically, consider the weight-2 newforms of level 11 and 17:

$$f_1(\tau) = e(\tau) - 2e(2\tau) - e(3\tau) + 2e(4\tau) + \mathcal{O}(e(5\tau)),$$
$$f_2(\tau) = e(\tau) - e(2\tau) - e(4\tau)^4 + \mathcal{O}(e(5\tau)).$$

The 19-th Hecke eigenvalue of $f_1$ and the 29-th one of $f_2$ vanish modulo 3. We conclude that $f_1 + f_2$ has a Ramanujan-type congruence modulo 3 on the arithmetic progression $19{\cdot}29\,(19{\cdot}29\,\mathbb{Z}+1)$, but neither on $19(19\mathbb{Z}+1)$ nor on $29(29\mathbb{Z}+1)$. By applying the Shintani lift one obtains a corresponding example in weight $\frac{3}{2}$.

**Example 3.18.** We next consider generalized Hecke eigenforms modulo $\ell$. For simplicity, assume that $\ell$ be a principal ideal. We identify $\ell$ with one of its ideal generators. Consider two newforms $f_1 \ne f_2$ of the same weight with $f_1 \equiv f_2 \not\equiv 0 \pmod{\ell}$. Then $g := (f_1 - f_2)/\ell$, if not divisible by $\ell$, yields a generalized Hecke eigenform for all primes not dividing the level and for which the Hecke eigenvalues of $f_1$ and $f_2$ do not satisfy a higher congruence. The following calculations can be verified, for instance using Sage [22].

We first discuss the case of integral weights in detail. Recall the newforms $f_1$ and $f_2$ of weight 2 and level 11 and 33 from Example 3.1. We set $\ell = 3$ and $g = (f_1 - f_2)/3$. We wish to make Lemma 3.9 explicit for $M = 35$. There is only one tuple of eigenvalues $\lambda = (\lambda_p)_{p\,|\,M}$ to consider, so we suppress it from our discussion. We order the divisors of $M$, and write $t = (t_5, t_7)$. The eigenvalues of $f_1$ and $f_2$ at 5 and 7 equal $1 \equiv -2 \pmod{3}$ and $-2 \equiv 4 \pmod{3}$. From this we see that

$$g_{(1,0)} \equiv \frac{1-(-2)}{3} f_1 = f_1 \pmod{3}, \quad g_{(0,1)} \equiv \frac{-2-4}{3} f_1 = -2 f_1 \pmod{3}.$$

We have $g_{(1,1)} \equiv 0 \pmod{3}$, so it does not contribute. Summarizing, we find that that $L(g;s)$ is congruent modulo 3 to

$$\frac{L^{35}(g,s)}{L_5(1,5^{-s})\,L_7(-2,7^{-s})} + \frac{L^{35}(f_1,s)\,5^{-s}}{L_5(1,5^{-s})^2\,L_7(-2,7^{-s})} + \frac{-2L^{35}(f_1,s)\,7^{-s}}{L_5(1,5^{-s})\,L_7(-2,7^{-s})^2}.$$

To check for Ramanujan-congruences, we have to inspect the coefficients of two Dirichlet series:

$$\frac{1}{L_5(1,5^{-s})\,L_7(-2,7^{-s})} \quad \text{and} \quad \frac{5^{-s}}{L_5(1,5^{-s})^2\,L_7(-2,7^{-s})} + \frac{-2\cdot 7^{-s}}{L_5(1,5^{-s})\,L_7(-2,7^{-s})^2}.$$

A calculation yields the following initial expansions modulo 3:

$$\begin{aligned} L_5(1,X)^{-1} &\equiv 1+X+2X^2+2X^4+2X^5+X^6+X^8+X^9+2X^{10}+\mathcal{O}\big(X^{12}\big),\\ L_7(-2,X)^{-1} &\equiv 1+X+2X^3+2X^4+X^6+X^7+2X^9+\mathcal{O}\big(X^{12}\big),\\ X\,L_5(1,X)^{-2} &\equiv X+2X^2+2X^3+X^4+2X^5+2X^6+2X^7+X^8+X^9+2X^{10}+\mathcal{O}\big(X^{12}\big),\\ X\,L_7(-2,X)^{-2} &\equiv X+2X^2+X^3+X^4+2X^5+X^6+2X^{10}+X^{11}+\mathcal{O}\big(X^{12}\big). \end{aligned}$$

We conclude that we have a Ramanujan-type congruence with gap

$$\forall n\in\mathbb{Z}, \gcd(n,35)=1 : c(g;\, 5^3 7^2 n) \equiv 0 \pmod{3}.$$

On the other hand, we find from a direct calculation that

$$c(g;\, 5^3) = -7 \not\equiv 0 \pmod{3}, \quad c(g;\, 7^2) = -4 \not\equiv 0 \pmod{3}.$$

To complete the picture, we observe that $g$ also has the following Ramanujan-type congruences with gap:

$$\forall n\in\mathbb{Z}\setminus 5\mathbb{Z} : c(g;\, 5^{11} n) \equiv 0 \pmod{3}.$$

In particular, there are Ramanujan-type congruences on $M\mathbb{Z}+\beta$, where $M$ is a prime power, but as opposed to the case of Hecke eigenforms modulo $\ell$ they do not give rise to all Ramanujan-type congruences.

**Example 3.19.** The generalized Hecke eigenform modulo 3 in weight 2 and level 33 does not give rise to a proper generalized Hecke eigenform in weight $\frac{3}{2}$ as explained in Example 3.1. Instead, Maeda's congruence between modular forms [12] can be used to also obtain an example in half-integral weight. The calculations in this example were performed using the package Hecke of the computer algebra system Oscar [7].

Observe that Maeda does not work with the Kohnen-plus space and that $\ell$ is not principal. We refer to Maeda's modular forms $F(-87)$ as $\tilde{f}_1$, $F((87+\sqrt{2304})/2$ as $\tilde{f}_2$, $f(-87)$ as $f_1$, and $f((87+\sqrt{2304})/2$ as $f_2$. We use the number field $K=\mathbb{Q}(\sqrt{2304})$ and the prime ideal $\ell=\mathcal{O}_{K,\ell}\{433, 172-\sqrt{2304}\}$. The newforms $\tilde{f}_1$ and $\tilde{f}_2$ have weight 8 and level 26. Recall the Hecke eigenvalues given by Maeda:

$$\begin{aligned} &c(\tilde{f}_1;2)=8,\ c(\tilde{f}_1;3)=-87, && c(\tilde{f}_1;5)=321, && c(\tilde{f}_1;7)=-181;\\ &c(\tilde{f}_2;2)=8,\ c(\tilde{f}_2;3)=\frac{87+\sqrt{2304}}{2}, && c(\tilde{f}_2;5)=\frac{215+5\sqrt{2304}}{2}, && c(\tilde{f}_2;7)=\frac{705-49\sqrt{2304}}{2}. \end{aligned}$$

The newforms $f_1$ and $f_2$ have weight $k = \frac{9}{2}$ and level 52 and transform like the theta-character. Their initial Fourier expansions are

$$\begin{aligned} f_1(\tau) &= 13e(2\tau) \;\; + 76e(5\tau) \;\; - 29e(6\tau) \;\; + \mathcal{O}\big(e(7\tau)\big), \\ f_2(\tau) &= 240e(2\tau) + (-650 + 10\sqrt{2305})e(5\tau) + (-20 + 52\sqrt{2305})e(6\tau) + \mathcal{O}\big(e(7\tau)\big). \end{aligned}$$

Maeda showed that $240 f_1 \equiv 13 f_2 \pmod{\ell}$. We verify that $\ell\, \frac{83+3\sqrt{2305}}{866} \subset \mathcal{O}_{K,\ell}$, and thus obtain a generalized eigenform with integral coefficients

$$\begin{aligned} g :&= \frac{83 + 3\sqrt{2305}}{866}\,(240 f_1 - 13 f_2) \\ &= (1520 + 80\sqrt{2305})e(5\tau) + (-6040 - 88\sqrt{2305})e(7\tau) + \mathcal{O}\big(e(7\tau)\big). \end{aligned}$$

There is an analogue of Lemma 3.9 that holds for modular forms that do not lie in the Kohnen-plus space. This analogue can be derived from Shimura's original relation of L-series, which encompasses (1.6). Instead of a fundamental discriminant $D$, we will need a nonzero square-free integer $T$ of sign $(-1)^{k-\frac{1}{2}}$. Since $k = \frac{9}{2}$ in Maeda's example, $T$ will be positive. We need to complete one more calculation, before we can provide an explicit description of the L-series of $g$. We adopt the notation in Lemma 3.9, write $t = (t_3, t_5, t_7)$, and compute the following congruences modulo $\ell$:

$$\begin{aligned} g_{(1,0,0)} &\equiv \frac{83 + 3\sqrt{2305}}{866}\, 240 \Big( -87 - \frac{87 + \sqrt{2305}}{2} \Big) f_1 &&\equiv (370 + 313\sqrt{2305})\, f_1, \\ g_{(0,1,0)} &\equiv \frac{83 + 3\sqrt{2305}}{866}\, 240 \Big( 321 - \frac{215 + 5\sqrt{2305}}{2} \Big) f_1 &&\equiv (120 + 120\sqrt{2305})\, f_1, \\ g_{(0,0,1)} &\equiv \frac{83 + 3\sqrt{2305}}{866}\, 240 \Big( -181 - \frac{705 - 49\sqrt{2305}}{2} \Big) f_1 &&\equiv (40 + 120\sqrt{2305})\, f_1. \end{aligned}$$

All higher $g_t$ vanish modulo $\ell$, since $f_1$ is an eigenform.

We now find that there are congruences modulo $\ell$ of the following type:

$$\begin{aligned} &\sum_{n=1}^{\infty} c(g;\, T n^2) n^{-s} \\ &\equiv c(g;\, T)\, \big(1 - \chi_T(5) 5^{3-s}\big)\big(1 - \chi_T(7) 7^{3-s}\big) \Bigg( \frac{\mathrm{L}^{35}\big(\mathrm{Sh}_T(g), s\big)}{\mathrm{L}_5(321, 5^{-s})\, \mathrm{L}_7(-181, 7^{-s})} \\ &\qquad + \frac{(120 + 120\sqrt{2305})\, \mathrm{L}^{35}\big(\mathrm{Sh}_T(f_1), s\big)\, 5^{-s}}{\mathrm{L}_5(321, 5^{-s})^2\, \mathrm{L}_7(-181, 7^{-s})} + \frac{(40 + 120\sqrt{2305})\, \mathrm{L}^{35}\big(\mathrm{Sh}_T(f_1), s\big)\, 7^{-s}}{\mathrm{L}_5(321, 5^{-s})\, \mathrm{L}_7(-181, 7^{-s})^2} \Bigg). \end{aligned}$$

Further formulas involving any combination of the primes 3, 5, and 7 follow the same pattern, but in light of the length of the resulting expressions we omit them.

Compared to the situation in the integral weight case, we now have to evaluate a much larger number of Dirichlet series, one for each value of $\chi_T(p) \in \{0, \pm 1\}$. Adding to this, the expansion length required to find Ramanujan-type congruences is too large to realistically give all coefficients of the relevant Dirichlet series. Instead, in what follows we give the exponents less than 1000 of $X$ for which the corresponding coefficients vanish modulo $\ell$.

At the prime 3, we find vanishing coefficients for exponents:

$$\frac{1}{\mathrm{L}_3(-87,X)} : 216, 433, 650, 867; \quad \frac{1-3^3X}{\mathrm{L}_3(-87,X)} : \text{none}; \quad \frac{1+3^3X}{\mathrm{L}_3(-87,X)} : \text{none};$$
$$\frac{X}{\mathrm{L}_3(-87,X)^2} : 0, 306, 329; \quad \frac{X(1-3^3X)}{\mathrm{L}_3(-87,X)^2} : 0, 494, 709, 817, 959; \quad \frac{X(1+3^3X)}{\mathrm{L}_3(-87,X)^2} : 0, 311, 554.$$

At the prime 5, we find vanishing coefficients for exponents:

$$\frac{1}{\mathrm{L}_5(321,X)} : 433, 867; \quad \frac{1-5^3X}{\mathrm{L}_5(321,X)} : 139, 573; \quad \frac{1+5^3X}{\mathrm{L}_5(321,X)} : 324, 758;$$
$$\frac{X}{\mathrm{L}_5(321,X)^2} : 0, 733; \quad \frac{X(1-5^3X)}{\mathrm{L}_5(321,X)^2} : 0, 454; \quad \frac{X(1+5^3X)}{\mathrm{L}_5(321,X)^2} : 0, 566, 670, 722, 843.$$

At the prime 7, we find vanishing coefficients for exponents:

$$\frac{1}{\mathrm{L}_7(-181,X)} : 433, 867; \quad \frac{1-7^3X}{\mathrm{L}_7(-181,X)} : 206, 640; \quad \frac{1+7^3X}{\mathrm{L}_7(-181,X)} : 60, 494, 928;$$
$$\frac{X}{\mathrm{L}_7(-181,X)^2} : 0, 444, 850; \quad \frac{X(1-7^3X)}{\mathrm{L}_7(-181,X)^2} : 0, 142, 661; \quad \frac{X(1+7^3X)}{\mathrm{L}_7(-181,X)^2} : 0, 100, 969.$$

At the end of this example, we show that the two Dirichlet series at the prime 3 whose coefficients of $X^m$, $m < 1000$, never vanish modulo $\ell$ in fact have only nonzero coefficient modulo $\ell$.

From the list of vanishing coefficients modulo $\ell$, we can extract the existence and non-existence of many Ramanujan-type congruences. For example, if $3 \nmid T$, then there is no Ramanujan-type congruences for $f_1$, $f_2$, or $g$ modulo $\ell$ on any arithmetic progression $T3^{2m}(3\mathbb{Z} + \beta_0)$ for any positive integer $\mathbb{Z}$. This follows since none of the coefficients of $X^m$ in the power series $(1 \pm 3^3X)\,\mathrm{L}_3(-87,X)^{-1}$ is zero modulo $\ell$; The cases $m \geq 1000$ will be discussed at the end of this example. Even in the case that $3 \mid T$, we find Ramanujan-type congruences for $f_1$ and $f_2$ on, for instance, $T3^{2\cdot 216}(3\mathbb{Z}+1)$, but $g$ does not have any such Ramanujan-type congruence, since none of the coefficients of $X^m$, $m < 1000$ in $\mathrm{L}_3(-87,X)^{-1}$ and $X\,\mathrm{L}_3(-87,X)^{-2}$ simultaneously vanish modulo $\ell$.

Next, we inspect the case of $T = 5$. We record that $c(g;T) \not\equiv 0 \pmod{\ell}$, $\chi_T(5) = 0$, and $\chi_T(7) = -1$. We conclude that we have the Ramanujan-type congruence

$$\forall n \in \mathbb{Z} : c\big(g; 5^{867} 7^{120}(35n+1)\big) \equiv 0 \pmod{\ell}.$$

We next want to exclude Ramanujan-type congruences on the arithmetic sup-progressions $5^{867}(5\mathbb{Z} + 7^{120}) = 5^{867}(5\mathbb{Z}+1)$ and $7^{120}(7\mathbb{Z} + 5^{867}) = 7^{120}(7\mathbb{Z}+6)$. To rule out the first one, we check that $c(g;5) \not\equiv 0 \pmod{\ell}$ and then use our previous list of vanishing coefficients of $X\,\mathrm{L}_5(321, X)^{-2} \pmod{\ell}$. The 433-th coefficient of this series does not vanish, disproving the first Ramanujan-type congruence. For the second one, we check that $c(g;6) \not\equiv 0 \pmod{\ell}$. Since $\chi_6(7) = -1$, we have to inspect the nonzero coefficients of $X(1+7^3X)\,\mathrm{L}_7(-181, X)^{-2} \pmod{\ell}$. The 60-th coefficient does not vanish, disproving the second Ramanujan-type congruence.

We finish this example by proving that neither $g$ nor $f_1$ or $f_2$ have Ramanujan-type congruence on any arithmetic progression of the form $T3^{2m}(3\mathbb{Z}+\beta)$, where $3 \nmid T$, $3 \nmid \beta$ except if we have one on $T(3\mathbb{Z}+\beta)$. This is equivalent to the statement that none of the coefficients of $(1 \mp 3^3 X)\,\mathrm{L}_3(-87, X)$ vanishes modulo $\ell$. We employ Proposition 3.12. To this end, we need to factor the L-polynomial $\mathrm{L}_3(-87,X) \pmod{\ell}$, which requires us to pass to the quadratic extension of $\mathrm{F}_{K,\ell} = \mathbb{F}_{433}$. Choosing a suitable generator $r$ for this extension $F = \mathrm{F}_{K,\ell}[r]$ with $r^2 + 432r + 5 = 0$, we find that

$$\mathrm{L}_3(-87, X) \;\equiv\; (1-\alpha_3 X)(1-\beta_3 X) \pmod{\ell}, \quad \alpha_3 = 126 + 94r,\, \beta_3 = 220 + 339r.$$

To apply the criterion in (3.19), we need the multiplicative orders:

$$\mathrm{ord}\Big(\frac{\alpha_p}{\beta_p}\Big) = 217, \quad \mathrm{ord}\Big(\frac{\beta_p - 3^3}{\alpha_p - 3^3}\Big) = 62,\, \mathrm{ord}\Big(\frac{\beta_p + 3^3}{\alpha_p + 3^3}\Big) = 434.$$

Since $62, 434 \nmid 217$, this disproves Ramanujan-type congruences modulo $\ell$ on arithmetic progressions $T3^{2m}(3\mathbb{Z}+\beta)$, where $3 \nmid T$, $3 \nmid \beta$, as desired.

# 4 Proofs of the Main Statements

*Proof of Theorem A.* This is a special case of Proposition 3.11 for $K = \mathbb{Q}$. ∎

The proof of Theorem B is provided in a separate subsection.

*Proof of Theorem C.* Statements (1) and (2) are special cases of Proposition 2.1. Statement (3) follows from Proposition 3.5. If $M_p = p^2$, Statement (4) is a special case of Proposition 2.3 when applied to $\mathrm{V}_{24} f$. If $M_p$ is at least a fourth power, then we employ Proposition 3.15 instead. Finally, Statement (5) is a special case of Proposition 3.16. ∎

*Proof of Theorem D.* Proposition 2.1 confirms the statement if $\beta'_0 \square \beta_0 \pmod{p}$.

If $f$ has a Ramanujan-type congruence on $p^2(p\mathbb{Z} + p) = p^3\mathbb{Z}$, then Proposition 3.5 implies that $f \equiv 0 \pmod{\ell}$, contradicting the assumption that the given Ramanujan-type congruence is maximal.

Since the Ramanujan-type congruence is maximal, Proposition 3.7 ensures that there is a fundamental discriminant in $\mathrm{D}_k(p^2(p\mathbb{Z}+\beta_0))$ such that $\mathrm{Sh}_D(f) \not\equiv 0 \pmod{\ell}$. Proposition 3.11 allows us to assume that $f$ is a Hecke eigenform modulo $\ell$ with respect to $\mathrm{T}^{\mathrm{cl}}_p$. This allows us to invoke Corollary 3.14 to finish the proof. ∎

*Proof of Theorem E.* The maximal Ramanujan-type congruences modulo $\ell$ for the partition function corresponds to a maximal Ramanujan-type congruences for $\eta^{-1}$ on $\ell M\mathbb{Z}+\beta-1$. Note that we can remove the factor 24 from $24\ell M\mathbb{Z}$, since the Fourier coefficients of $\eta^{-1}$ are supported on $\frac{-1}{24} + \mathbb{Z}$. We thus obtain a maximal Ramanujan-type congruences for $f_{\ell,\delta}$, defined in the introduction after the statement of Theorem E. If $\delta = 0$ this Ramanujan-type congruences lives on $M\mathbb{Z}+(\beta-1)/\ell$, and if $\delta = -1$, we have such a congruence on $\ell M\mathbb{Z} + \beta - 1$. Using Proposition 3.11 we obtain the same Ramanujan-type congruence (not necessarily maximal) for a nontrivial Hecke eigenform modulo $\ell$, which we denote by $f \not\equiv 0 \pmod{\ell}$. Observe that $f$ has the same weight and character as $f_{\ell,\delta}$. Lemma 3.6 shows that $\mathrm{Sh}_D(f) \equiv 0 \pmod{\ell}$ for all $D \in \mathrm{D}_k(M\mathbb{Z} + (\beta-1)/\ell)$ if $\delta = 0$ and $D \in \mathrm{D}_k(\ell M\mathbb{Z} + \beta - 1)$ if $\delta = -1$. Now the same induction argument as given later in the proof of Theorem B, shows the desired congruences for $f$. ∎

*Proof of Theorem F.* As in the proof of Theorem E, we obtain a Hecke eigenform modulo $\ell$ that has a Ramanujan-type congruence on $q^3\mathbb{Z} + (\beta-1)/\ell$ or $q^4\mathbb{Z} + (\beta-1)/\ell$, respectively, if $\delta = 0$, and on $\ell q^3\mathbb{Z} + \beta - 1$ or $\ell q^4\mathbb{Z} + \beta - 1$ if $\delta = -1$. Write $\lambda_p$ for its Hecke eigenvalues modulo $\ell$.

We also find that $f_{\ell,\delta}$ has a maximal Ramanujan-type congruence on the respective arithmetic progression. By Proposition 2.1, we have $q^2 \mid (\beta-1)$, which implies that $\beta_0$ is integral. Moreover, we have $q \nmid (\beta_0 - 1)$ by Proposition 3.5 in the first case of the proposition, and $q \mid (\beta_0 - 1)$ in the second one again by Proposition 2.1.

In the first case of the proposition, Proposition 3.12 with $\epsilon = \pm 1$ and $m = 1$ asserts that we have the relation $\lambda_p \equiv \pm p^{k-\frac{3}{2}} \pmod{\ell}$. We have have $k - \frac{3}{2} = \frac{\ell-5}{2}$ if $\delta = 0$, and $k - \frac{3}{2} \equiv -2 \pmod{\ell - 1}$ if $\delta = -1$. If $\delta = -1$, we obtain the desired congruences directly, and if $\delta = 0$, we calculate that the square of $\lambda_p$ is congruent to $p^{-2}$ modulo $\ell$.

In the second case, Proposition 3.12 with $\epsilon = 0$ and $m = 1$ directly yields the desired congruence $\lambda_p \equiv 0 \pmod{\ell}$, using that $m \not\equiv -1 \pmod{\ell}$. ∎

**4.1 Proof of Theorem B** Since the Ramanujan-type congruence of $f$ is maximal, the discussion at the beginning of Section 3 shows that we have $\gcd(M,6) = 1$. We therefore can replace $f$ by $\mathrm{V}_{24}\, f$ and work with the theta-character instead of the eta-character. For simplicity, we replace $N$ by $\mathrm{lcm}(N,4)$ and assume that $4 \,|\, N$. Since $M$ is odd, we have $M \,|\, M_{\mathrm{sf}} \beta$ by Proposition 2.1, and hence $\beta_1 \in M_{\mathrm{fd}} \mathbb{Z}$.

Assume that $\mathrm{Sh}_D(f) \not\equiv 0 \pmod{\ell}$ for some $D \in \mathrm{D}_k(M\mathbb{Z} + \beta)$. We want to show that the eigenvalues $\lambda_p$ modulo $\ell$ of $f$ satisfy the congruences in Proposition 3.11. Given $p \,|\, M$ set

$$f_p \;:=\; \sum_{\substack{n \in \mathbb{Z} \\ n \,\square\, \beta \pmod{M_p^{\#}}}} c(f;\, n) e(n\tau),$$

which has a Ramanujan-type congruence on $M_p \mathbb{Z} + \beta$ and by Lemma 2.2 is a modular form for the character $\chi \chi_\theta^r$ of $\widetilde{\Gamma}_0(M_p^{\#2} N)$. We have $\mathrm{Sh}_D(f_p) \not\equiv 0 \pmod{\ell}$ for all $p \,|\, M$.

For any nonnegative integer $t$, we set

$$f_{p,t} \;:=\; f_p \big| \big(\mathrm{T}_p^{\mathrm{cl}} - \lambda_p\big)^t.$$

By Proposition 3.11, $f_{p,t}$ has a Ramanujan-type congruence module $\ell$ on $M_p \mathbb{Z} + \beta$ for every $t$. Further, we define $\mathrm{Dsupp}_{p,t}$ as the set

$$\Big\{ D \text{ fund. disc.} : \mathrm{sgn}(D) = (-1)^{k - \frac{1}{2}},\, |D| \,\square\, \beta_1 \pmod{M_p},\, \mathrm{Sh}_D(f_{p,t}) \not\equiv 0 \pmod{\ell} \Big\}.$$

For given $p$, we have $\mathrm{Dsupp}_{p,0} \neq \emptyset$ by our assumption that $\mathrm{Sh}_D(f_p) \not\equiv 0 \pmod{\ell}$. For sufficiently large $t$, we have $\mathrm{Dsupp}_{p,t} = \emptyset$, since then $f_{p,t} \equiv 0 \pmod{\ell}$. In particular, there is a nonnegative integer $d$ such that we have $\mathrm{Dsupp}_{p,d+1} = \emptyset$, while $\mathrm{Dsupp}_{p,d}$ is nonempty. We conclude that $f_{p,d}$ satisfies the assumptions of Proposition 3.12. This finishes the proof if there is $D \in \mathrm{D}_k(M\mathbb{Z} + \beta)$ with $\mathrm{Sh}_D(f) \not\equiv 0 \pmod{\ell}$.

Assume that we have $\mathrm{Sh}_D(f) \equiv 0 \pmod{\ell}$ for all $D \in \mathrm{D}_k(M\mathbb{Z} + \beta)$. We have to show that $\mathrm{Sh}_{\tilde{D}}(f) \equiv 0 \pmod{\ell}$ for all $\tilde{D} \in \widetilde{D}_k(M\mathbb{Z} + \beta)$. If $\ell \,|\, 2$, we have $\widetilde{D}_k(M\mathbb{Z} + \beta) = \mathrm{D}_k(M\mathbb{Z} + \beta)$ by Lemma 3.6. In this case, the proof is complete. We assume in the remainder that $\ell \nmid 2$. Fix any $\tilde{D} \in \widetilde{D}_k(M\mathbb{Z} + \beta)$. We can and will assume that $\tilde{D} \not\in \mathrm{D}_k(M\mathbb{Z} + \beta)$.

We are assuming that $\mathrm{Sh}_D(f) \equiv 0 \pmod{\ell}$ for all $D \in \mathrm{D}_k(M\mathbb{Z} + \beta)$ and Proposition 3.7 shows that we have the congruence

$$\forall n_0 \in \mathbb{Z}, n_0 \,\square\, \beta_0 \pmod{M_0} : c(f;\, M_{\mathrm{fd}} n_0) \equiv 0 \pmod{\ell}.$$

In particular, we find the Ramanujan-type congruence

$$\forall n_0 \in \mathbb{Z} : c\big(f;\, M_{\mathrm{fd}}(M_0 n_0 + \beta_0)\big) \equiv 0 \pmod{\ell}.$$

We record for later purpose that we have $\gcd(\beta_0, M_{\mathrm{fd}}) = 1$ by Proposition 3.5, since the Ramanujan-type congruence of $f$ is maximal.

Let $P$ be the set of primes $p$, $p \,|\, M$, for which $|\tilde{D}|$ and $M_{\mathrm{fd}}\beta_0$ are not in the same square-class modulo $p^2$. Set $M_P^\# = \prod M_p$ where the product runs over $p \,|\, M$ with $p \not\in P$. We write $B$ for the set of $\tilde{\beta}_0 \in \mathbb{Z}$ with $p \nmid \tilde{\beta}_0$ for all $p \in P$, $\tilde{\beta}_0 \equiv \beta_0 \pmod{\gcd(M_P^\#, M_0)}$, and

$$\forall n_0 \in \mathbb{Z} : c\big(f;\, M_{\mathrm{fd}}(M_0 n_0 + \tilde{\beta}_0)\big) \equiv 0 \pmod{\ell}.$$

We claim that given any choice of signs $\epsilon_p \in \{\pm 1\}$, $p \in P$, there is some $\tilde{\beta}_0 \in B$ satisfying the Jacobi symbol condition $(\tilde{\beta}_0 / p) = \epsilon_p$. We prove the statement by induction on the number of primes $p \in P$ with $\epsilon_p \ne (\beta_0 / p)$. The base case of the induction holds, since we have $\beta_0 \in B$ by the Ramanujan-type congruence that we have already discovered. Fix some choice of $\epsilon_p$, $p \in P$, and assume that there is a prime $q \in P$ with $\epsilon_q \ne (\beta_0 / q)$. By induction, there is $\tilde{\beta}'_0 \in B$ with $(\tilde{\beta}'_0 / p) = \epsilon_p$ for all $p \in P$, $p \ne q$, and $(\tilde{\beta}'_0 / q) = (\beta_0 / q) = -\epsilon_q$. We have a Ramanujan-type congruence modulo $\ell$ on $M_{\mathrm{fd}}(M_0 \mathbb{Z} + \tilde{\beta}'_0)$. Further, we have $\gcd(\tilde{\beta}'_0, M_{\mathrm{fd}}) = \gcd(\beta_0, M_{\mathrm{fd}}) = 1$ and by assumption we have $\ell \nmid 2$. Therefore can apply Proposition 2.3 to obtain the desired $\tilde{\beta}_0$.

Now fix $\tilde{\beta}_0 \in B$ with $M_{\mathrm{fd}}\tilde{\beta}_0 \,\square\, |\tilde{D}| \pmod{p^2}$ for all $P$. We now invoke Proposition 2.1 and find that

$$\forall n_0 \in \mathbb{Z},\, n_0 \,\square\, \tilde{\beta}_0 \pmod{M_0} : c(f;\, M_{\mathrm{fd}} n_0) \equiv 0 \pmod{\ell}.$$

Since the Ramanujan-type congruence of $f$ is maximal, $M_0 = M_{\mathrm{sf}}$. In particular, $f$ satisfies all Ramanujan-type congruences of cube-free period as given in Theorem B. Lemma 3.6 implies that $c(f; |\tilde{D}|) \equiv 0 \pmod{\ell}$, since $\tilde{D} \in \mathrm{D}_k(M_{\mathrm{fd}}(M_0 \mathbb{Z} + \tilde{\beta}_0))$, confirming the congruence for fundamental discriminants.

[1] S. Ahlgren, P. B. Allen, and S. Tang. *Congruences like Atkin's for the partition function.* arXiv:2203.11273. 2022.

[2] S. Ahlgren, O. Beckwith, and M. Raum. *Scarcity of congruences for the partition function.* arXiv:2006.07645. 2020.

[3] S. Ahlgren and K. Ono. "Congruence properties for the partition function". *Proc. Natl. Acad. Sci. USA* 98.23 (2001).

[4] A. O. L. Atkin. "Multiplicative congruence properties and density problems for $p(n)$". *Proc. London Math. Soc. (3)* 18 (1968).

[5] J. H. Bruinier. "Nonvanishing modulo $l$ of Fourier coefficients of half-integral weight modular forms". *Duke Math. J.* 98.3 (1999).

[6] N. Dummigan. *Lifting Congruences to Half-Integral Weight.* http://neil-dummigan.staff.shef.ac.uk/papers.html. 2021.

[7] C. Fieker, W. Hart, T. Hofmann, and F. Johansson. *Oscar/Hecke (Version 0.8.6).* https://oscar.computeralgebra.de/. 2021.

[8] F. Johansson. "Efficient implementation of the Hardy-Ramanujan-Rademacher formula". *LMS J. Comput. Math.* 15 (2012).

[9] L. J. P. Kilford. "Generating spaces of modular forms with $\eta$-quotients". *JP J. Algebra Number Theory Appl.* 8.2 (2007).

[10] W. Kohnen. "Fourier coefficients and modular forms of half-integral weight". *Math. Ann.* 271 (1985).

[11] H. Maass. *Lectures on modular functions of one complex variable.* Notes by Sunder Lal. Tata Institute of Fundamental Research Lectures on Mathematics, No. 29. Bombay: Tata Institute of Fundamental Research, 1964.

[12] Y. Maeda. "A congruence between modular forms of half-integral weight". *Hokkaido Math. J.* 12 (1983).

[13] S. Niwa. "Modular forms of half integral weight and the integral of certain theta-functions". *Nagoya Math. J.* 56 (1975).

[14] K. Ono. "Distribution of the partition function modulo $m$". *Ann. of Math. (2)* 151.1 (2000).

[15] C.-S. Radu. "A proof of Subbarao's conjecture". *J. Reine Angew. Math.* 672 (2012).

[16] C.-S. Radu. "Proof of a conjecture by Ahlgren and Ono on the non-existence of certain partition congruences". *Trans. Amer. Math. Soc.* 365.9 (2013).

[17] S. Ramanujan. "Congruence properties of partitions." *Proc. Lond. Math. Soc. (2)* 18 (1920).

[18] M. Raum. *Relations among Ramanujan-Type Congruences I.* arXiv:2010.06272v2. 2020.

[19] G. Shimura. "On modular forms of half integral weight". *Ann. of Math. (2)* 97 (1973).

[20] T. Shintani. "On construction of holomorphic cusp forms of half integral weight". *Nagoya Math. J.* 58 (1975).

[21] F. Strömberg. "Weil representations associated with finite quadratic modules". *Math. Z.* 275.1-2 (2013).

[22] The Sage Developers. *SageMath, the Sage Mathematics Software System (Version 9.2).* https://www.sagemath.org. 2020.

[23] J.-L. Waldspurger. "Sur les coefficients de Fourier des formes modulaires de poids demi-entier". *J. Math. Pures Appl. (9)* 60.4 (1981).

[24] R. L. Weaver. "New congruences for the partition function". *Ramanujan J.* 5.1 (2001).

[25] Y. Yang. "Modular forms of half-integral weights on $\mathrm{SL}(2,\mathbb{Z})$". *Nagoya Math. J.* 215 (2014).

Chalmers tekniska högskola och Göteborgs Universitet, Institutionen för Matematiska vetenskaper, SE-412 96 Göteborg, Sweden
E-mail: martin@raum-brothers.eu
Homepage: http://raum-brothers.eu/martin